\documentclass[aos,MSNbibl,seceqn,nameyear,dvips]{arximspdf}

\doi{10.1214/13-AOS1152} 
\volume{42}
\issue{1}
\pubyear{2014}
\firstpage{1}
\lastpage{28}

\makeatletter
\newcommand{\rrVert}{\Vert}
\newcommand{\rrvert}{\vert}
\newcommand{\llVert}{\Vert}
\newcommand{\llvert}{\vert}
\newcommand{\widebar}{\overline}
\newproclaim{assumption}{Assumption}
\newtheorem{theorem}{Theorem}
\newtheorem{proposition}{Proposition}
\newproclaim{definition}{Definition}
\newtheorem{lemma}{Lemma}
\newproclaim{remark}{Remark}
\makeatother

\begin{document}
\begin{frontmatter}

\title{Adaptive estimation under single-index constraint in a regression model\thanksref{T1}}
\runtitle{Adaptation under single-index constraint}

\begin{aug}
\author[A]{\fnms{Oleg} \snm{Lepski}\ead[label=e1]{Oleg.Lepski@cmi.univ-mrs.fr}}
\and
\author[B]{\fnms{Nora} \snm{Serdyukova}\corref{}\ead[label=e3]{Nora.Serdyukova@gmail.com}}
\runauthor{O. Lepski and N. Serdyukova}
\affiliation{Aix-Marseille Universit\'e and Universidad de Concepci\'{o}n}
\address[A]{Laboratoire d'Analyse\\
Topologie, Probabilit\'es UMR 7353\\
Aix-Marseille Universit\'e\\
39, rue F. Joliot Curie \\
13453 Marseille\\
France\\
\printead{e1}}
\address[B]{Departamento de Estad\'{\i}stica\\
Facultad de Ciencias F\'{\i}sicas y Matem\'{a}ticas\\
Universidad de Concepci\'{o}n\\
Avda. Esteban Iturra s/n---Barrio Universitario\\
Concepci\'{o}n, Regi\'{o}n VIII\\
Chile\\
\printead{e3}}
\end{aug}
\thankstext{T1}{Supported in part by the DFG FOR 916 Statistical
Regularisation and Qualitative Constraints.}

\received{\smonth{11} \syear{2012}} \revised{\smonth{7} \syear{2013}}

%
\begin{abstract}
The problem of adaptive multivariate function estimation in the
single-index regression model with random design and weak assumptions
on the noise is investigated. 
A novel estimation procedure that adapts simultaneously to the unknown
index vector and the smoothness of the link function by selecting from
a family of specific kernel estimators is proposed. We establish a
pointwise oracle inequality which, in its turn, is used to judge the
quality of estimating the entire function (``global'' oracle
inequality). Both the results are applied to the problems of pointwise
and global adaptive estimation over a collection of H\"older and
Nikol'skii functional classes, respectively.
\end{abstract}

%
\begin{keyword}[class=AMS]
\kwd[Primary ]{62G05}
\kwd[; secondary ]{62G08}
\kwd{62G20}
\end{keyword}
\begin{keyword}
\kwd{Adaptive estimation}
\kwd{lower bounds}
\kwd{minimax rate}
\kwd{nonparametric regression}
\kwd{oracle inequalities}
\kwd{single-index model}
\kwd{structural adaptation}
\end{keyword}

\end{frontmatter}

\section{Introduction}\label{SectionIntro}

This paper deals with multivariate functions estimation. For the
proposed estimator we establish local as well as global oracle
inequalities and show how to use them for deriving minimax adaptive
results.\vspace*{9pt}

\textit{Model and setup.}
We observe $ (X_1, Y_1), \ldots, (X_n, Y_n) \in\mathbb{R}^d\times
\mathbb{R}$ following
%
\begin{equation}
\label{model-regression} Y_i = F(X_i) + \varepsilon_i,
\qquad i=1, \ldots, n,
\end{equation}
where $d\ge2$, the noise $ \{ \varepsilon_i\}_{i=1}^n $ are i.i.d.
centered random variables satisfying a tail probability condition
(Assumption \ref{assassumption-on-noise}), and the design points $ \{
X_i \}_{i=1}^n $ are independent random vectors with common density $ g
$ with respect to the Lebesgue measure. The sequences $ \{
\varepsilon_i\}_{i=1}^n $ and $ \{ X_i \}_{i=1}^n $ are assumed to be
independent. The density $g$ is known, however, in
Section~\ref{secunknowndensity} we discuss how to extend our results to
the case of unknown design density.

In addition, we assume that the function $ F\dvtx
\mathbb{R}^d\to\mathbb {R}$ has a single-index structure, that is,
there exist unknown $ f\dvtx  \mathbb{R}\to\mathbb{R} $ and $ \theta^*
\in\mathbb{R}^{d} $ such that
%
\begin{equation}
\label{single-index} F(x) = f\bigl(x^{\top} \theta^*\bigr).
\end{equation}
A minimal technical assumption about $ f $ is that it belongs to some
H\"older ball, yet the knowledge of this ball will not be required for
the proposed estimation procedure; see the discussion after Assumption
\ref{asstechnical} for more details.

The paper aims at estimating the entire function $ F $ on $ [-1/2,
1/2]^2 $ or its value $ F(t) $,  $t \in[-1/2, 1/2]^2 $, from the data $
\{ (X_i, Y_i)\}_{i=1}^n$ without any prior knowledge about the nuisance
parameters $ f(\cdot) $ and $ \theta^{*} $. The unit square is chosen
for notation convenience; and all the results remain true when $ [-1/2,
1/2]^2$ is replaced by an arbitrary bounded interval of $\mathbb{R}^2$.

Throughout the paper we adopt the following notation. The joint
distribution of the sequence $ \{ (X_i,Y_i)\}_{i=1}^n$ will be denoted
by $\mathbb{P}^{(n)}_{F}$, and those of $ \{
(X_i,\varepsilon_i)\}_{i=1}^n$ by~$\mathbb{P}^{(n)}_{X,\varepsilon}$.
In addition, $\mathbb{P}^{(n)}_{X}$ and
$\mathbb{P}^{(n)}_{\varepsilon}$ stand for the marginal distributions
of $ \{ X_i\}_{i=1}^n $ and $ \{ \varepsilon_i\}_{i=1}^n $,
respectively.

To judge the quality of estimation, we use either the risk determined
by the $ L_r  \operatorname{norm}$, $ \| \cdot\|_r $, on $ [-1/2, 1/2]^2
$ with $ r \in[1, \infty) $:
%
\begin{equation}
\label{globalrisk} \mathcal R_r^{(n)} (\widehat F, F) =
\mathbb{E}^{(n)}_{F} \| \widehat F - F \|_r,
\end{equation}
a ``global'' risk; or the ``pointwise'' risk defined as follows:
%
\begin{equation}
\label{localrisk} \mathcal R_{r,t}^{(n)} (\widehat F, F) = \bigl(
\mathbb{E}^{(n)}_{F} \bigl|\widehat F (t) - F(t)\bigr|^r
\bigr)^{1/r},\qquad t\in[-1/2,1/2]^2.
\end{equation}
Here $ \widehat F(\cdot) $ is an estimator, that is, an $ \{ (X_i, Y_i
)\}_{i=1}^n$-measurable function, and $ \mathbb{E}_F^{(n)} $ denotes
the mathematical expectation with respect to $ \mathbb{P}_F^{(n)} $.

All the results established in the paper, except the lower bound given
in Theorem~\ref{thpointwise-adaptation-lower}, are obtained for $ d=2
$. The principal difficulties with the case of arbitrary dimension are
commented in Remark \ref{remdifficulties}. It is noteworthy that the
single-index modeling, even if $d=2$, is a direct generalization of the
univariate regression model. Therefore, our results, mainly presented
in Section~\ref{subsecadaptive-estimation}, generalize in several
directions the existing ones obtained for the univariate random design
regression (see
the discussion after Theorem~\ref{thglobal-adaptation}).

\subsection*{Main assumptions} Let us formulate the principal
assumptions used in the sequel. They are imposed on the distributions
of the design and noise variables as well as on the approximation
property of the link function.

\begin{assumption}
\label{assassumption-on-noise} The random variable $\varepsilon_1$ has
a symmetric distribution with density $p$ with respect to the Lebesgue
measure. Moreover, there exist $\Upsilon>0$, $\Omega\in(0,1]$, and
$\omega>0$ such that
\[
p\in\mathfrak{P}= \biggl\{\ell\dvtx  \mathbb{R}\to\mathbb{R}_+ \bigg| \int
_{x}^{\infty}\ell(y)\,\mathrm{d}y\le\Upsilon
e^{-\Omega x^{\omega}}\ \forall x\ge0 \biggr\}.
\]
\end{assumption}

The assumption holds, for example, for the Gaussian, Laplace or, more
generally, for the symmetrized Weibull distribution. In the following,
the functional class $\mathfrak{P}$ is considered as fixed.

\begin{assumption}
\label{assassumption-on-design} There exists $ \underline{g}\in(0,1) $
such that $ \inf_{x \in[-3, 3]^2 }g(x) \ge\underline{g}$.
\end{assumption}

The assumption holds obviously if the design points are uniformly
distributed on any bounded Borel set containing $ [-3, 3]^2 $.
The imposed condition is ``fitted'' to the estimation over $[-1/2,
1/2]^2$ that explains the set $[-3, 3]^2$. When estimating over a
rectangle $[a,b]\times[c,e]\in\mathbb{R}^2$, the infimum should be
taken over $ [a-5/2,b+5/2]\times[c-5/2,e+5/2] $. If $ M $ from
Assumption~\ref{asstechnical} below is known, the above condition can
be relaxed to $ [a-2,b+2]\times[c-2,e+2] $. We also remark that
independently of the values $a,b,c,e$
Assumption~\ref{assassumption-on-design} is fulfilled if
$g\in\mathbb{C} (\mathbb{R}^2 )$ and $g(x)>0$ for any $x\in
\mathbb{R}^2$.

\begin{assumption}
\label{asstechnical} There exist $\beta_0\in(0,1)$ and $M>0$ such that
\[
f\in\mathbb{F}(\beta_0,M)= \biggl\{U\dvtx \mathbb{R}\to\mathbb{R} \bigg| \|U\|
_\infty+ \sup_{y_1,y_2\in\mathbb{R}}\frac{|U(y_1)-U(y_2)|}{|y_1-y_2|^{\beta
_0}}\le M \biggr\}.
\]
\end{assumption}

The latter assumption guarantees that the link function is smooth.
However, it is important to emphasize that $\beta_0$ and $M$ are not
supposed to be known {a priori}. In particular, they are not
involved in our estimation procedure. On the other hand, both the
parameters restrict the minimal sample size needed to justify the
theoretical results involved. Set for any $n\in\mathbb{N}^*$
%
\begin{equation}
\label{eqbandwidts} h_{\min}=n^{-1} \ln^{1+2/\omega}(n),\qquad
\mathfrak {h}=\sqrt{n^{-1} \ln^{1+1/\omega}(n)}.
\end{equation}
In the sequel it will be assumed that $n\ge n_0$, where
%
\begin{equation}
\label{eqrestriction-on-sample-size} n_0=\inf \bigl\{m\in\mathbb{N}^* | (M\vee1)\max \bigl
\{ \mathfrak{h}^{\beta_0}, \ln^{1/\omega}(n)h_{\min}^{\beta
_0}\bigr\}\le1\ \forall n\ge m \bigr\}.
\end{equation}
To finish this section, we remark that all the presented results remain
true if one assumes that $ f\in\mathbb{F}(0,M) $, that is, is uniformly
bounded, and $M$
is known. 

\subsection*{Objectives}
For clarity of presentation, it is
assumed that the index vector \mbox{$\theta^*\in\mathbb{S}^{d-1}$}, where
$\mathbb{S}^{d-1}$ stands for the unite sphere in $ \mathbb{R}^d$.
However, in Section~\ref{secextensions} it is shown that our results
can be extended to the case $\theta^*\in\mathbb{R}^2$.

The goal of our studies is at least threefold. We first seek an
estimation procedure $ \widehat{F}(t) $,  $t\in[-1/2,1/2]^{2} $, for $
F $ which could be applicable to any function $F$ satisfying
assumption~(\ref{single-index}). Moreover, we would like to bound the
risk of this estimator uniformly over the set
$\mathbb{F}(\beta_0,M)\times\mathbb{S}^{1}$. More precisely, we want to
establish for $ \widehat{F}(t) $ the so-called local oracle
inequality---at any point $ t\in[-1/2,1/2]^{2} $ the risk of $
\widehat{F}(t) $ should be bounded as follows:
%
\begin{equation}
\label{eqlocal-oracle-intro} \mathcal R_{r,t}^{(n)} (\widehat F, F)\le
C_r A^{(n)}_{f,\theta
^*}(t)\qquad\forall f\in\mathbb{F}(
\beta_0,M),\ \forall\theta ^*\in\mathbb{S}^{1}.
\end{equation}
Here $A^{(n)}_{f,\theta^*}(\cdot)$ is completely determined by the
function $f$, vector~$\theta^*$ and observations number $ n $, while
$C_r$ is a constant independent of $F$ and $n$.\vadjust{\goodbreak}

After being established, the local oracle inequality allows deriving
minimax adaptive results for the function estimation at a given point.
Indeed, let $ \{ \mathbb{F}(\gamma)$, $\gamma\in\Gamma\} $ be a
collection of functional classes such that
$\bigcup_{\gamma\in\Gamma}\mathbb{F}(\gamma)\subseteq\mathbb
{F}(\beta_0,M)$. For any $\gamma\in\Gamma$ define
\[
\phi_{n} (\gamma)=\inf_{\widetilde{F}}\sup
_{(f,\theta^*)\in\mathbb
{F}(\gamma
)\times\mathbb{S}^{1}}\mathcal R_{r,t}^{(n)} (\widetilde F,
F),
\]
where the infimum is taken over all possible estimators. The quantity
$\phi_n(\gamma)$ is the minimax risk on $\mathbb{F}(\gamma)\times
\mathbb{S}^{1}$. In the framework of minimax adaptive estimation, the
task is to construct an estimator $F^*$ such that for any
\mbox{$\gamma\in\Gamma$}
%
\begin{equation}
\label{eqlocal-adaptation-intro} \sup_{(f,\theta^*)\in\mathbb{F}(\gamma)\times\mathbb
{S}^{1}}\mathcal R_{r,t}^{(n)}
\bigl(F^*, F\bigr) \asymp \phi_{n}(\gamma),\qquad n\to\infty.
\end{equation}
The estimator $F^*$ satisfying (\ref{eqlocal-adaptation-intro}) is
called optimally rate adaptive over the collection $
\{\mathbb{F}(\gamma), \gamma\in\Gamma\} $. Subsequently, let
(\ref{eqlocal-oracle-intro}) be proved; and let for
any~\mbox{$\gamma\in\Gamma$}
\[
\sup_{(f,\theta^*)\in\mathbb{F}(\gamma)\times\mathbb{S}^{1}
}A^{(n)}_{f,\theta^*}(t)\asymp
\phi_{n}(\gamma),\qquad n\to\infty.
\]
Then one can assert that the estimator $\widehat{F}$ is adaptive over
$\{\mathbb{F}(\gamma), \gamma\in\Gamma\}$.

Thus, the first step is to prove (\ref{eqlocal-oracle-intro}). To the
best of our knowledge, such results do not exist in the context of
regression with random design not only under the single-index
constraint, but also in univariate regression.

Next, (\ref{eqlocal-oracle-intro}) is applied to minimax adaptive
estimation over H\"older classes,
$\{\mathbb{F}(\gamma)=\mathbb{H}(\beta,L),  \gamma=(\beta,L) \} $; see
Section~\ref{subsecadaptive-estimation} for pertinent definitions. We
will find the minimax rate over $\mathbb{H}(\beta,L)\times\mathbb
{S}^{1}$ and prove that $\widehat{F}$ achieves it, that is, is
optimally rate adaptive. This result is quite surprising because, if
$\theta^*$ is fixed, say, $\theta^*=(1,0)^{\top}$, it is well known
that an optimally adaptive estimator does not exist; see
\citet{Lep1990} for the Gaussian white noise model,
\citet{BrownLow} for density estimation, and \citet{Gaiffas}
for regression.

Local oracle inequality (\ref{eqlocal-oracle-intro}) allows us to bound
from above the ``global'' risk as well. Indeed, for any $ r \ge1 $, in
view of Jensen's inequality and Fubini's theorem, $   [ \mathcal
R_r^{(n)} (\widehat F, F) ]^r\le\mathbb{E}^{(n)}_{F}  \| \widehat F - F
\|^{r}_r= \|\mathcal R_{r,\cdot}^{(n)} (\widehat F, F) \| _r^r $ and,
therefore,
%
\begin{equation}
\label{eqglobal-oracle-intro} \mathcal R_r^{(n)} (\widehat F, F)\le
C_r \bigl\|A^{(n)}_{f,\theta
^*} \bigr\|_r.
\end{equation}
Inequality (\ref{eqglobal-oracle-intro}) is called the global oracle
inequality, and in the considered framework it supplies new results. As
local oracle inequality (\ref{eqlocal-oracle-intro}) is a powerful tool
for deriving minimax adaptive results in pointwise estimation, so
inequality (\ref{eqglobal-oracle-intro}) can be used for constructing
adaptive
estimators of $ F $.

We will consider a collection of Nikol'skii classes $\mathbb
{N}_p(\beta,L)$ (see Definition~\ref{defnikolskii-class}), where
$\beta,L>0$ and $1\le p<\infty$. When considering these classes, we aim
at estimating functions with inhomogeneous smoothness. This means that
the underlying function can be very regular on some parts of its domain
and rather irregular on the other sets. We will compute bounds for
\[
\sup_{(f,\theta^*)\in\mathbb{N}_p(\beta,L)\times\mathbb
{S}^{1}} \bigl\| A^{(n)}_{f,\theta^*}
\bigr\|_r
\]
and show that, if $(2\beta+1)p<r$, the rate of convergence is the
minimax rate over $\mathbb{N}_p(\beta,L)\times\mathbb{S}^{1}$. This
means that our estimator $\widehat{F}$ is optimally rate adaptive over
the collection $ \{\mathbb{N}_p(\beta,L)\times\mathbb{S}^{1},
\beta>0,L>0 \}$ whenever $(2\beta+1)p<r$. In the case $(2\beta+1)p\ge
r$, we will show that the latter bound differs from the bound on the
minimax risk by a logarithmic factor. Following the contemporary
language, we say that $\widehat{F}$ is ``nearly'' adaptive. The
construction of an optimally rate adaptive over the entire range of the
Nikol'skii classes estimator under the single-index constraint~(\ref{single-index}) is an open question.

All presented results are completely new. The adaptive estimation under
the $L_r$ loss and single-index constraint, except the case $r=2$ in
\citet{Lecue}, was not studied. Note, however, that the cited
result was obtained under the Gaussian errors model and over the
H\"older classes that do not admit the consideration of functions with
inhomogeneous smoothness.

\subsection*{Remarks}
It turns out that the adaptation to the unknown $ \theta^* $ and $
f(\cdot) $ can be viewed as selecting from a special family of kernel
estimators in the spirit of that of \citet{Lep1990},
\citet{KLP2001}, \citet{GoldLep2008}. However, our selection
rule is quite different from the aforementioned proposals, and it
allows us to solve the problem of minimax adaptive estimation under the
$ L_r $ losses over a collection of Nikol'skii classes.

It is worth mentioning that the single-index model is particularly
popular in econometrics [see, e.g., \citet{Horowitz1998},
\citet{MADDALA}]. The estimation, nevertheless, is usually
performed under smoothness assumptions on the link function. One
usually uses the $ L_2 $ losses, and the available methodology is based
on these restrictions. To the best of our knowledge, the only
exceptions are \citet{Golubev1992} for the minimax estimation
under the projection pursuit constraints, and \citet{GoldLep2009}
for adaptation to unknown smoothness and structure.

\subsection*{Organization of the paper} In Section~\ref{subsecoracle-approach} we present our selection rule and establish
for it local and global oracle inequalities. Section~\ref{subsecadaptive-estimation} is devoted to the application of these
results to minimax adaptive estimation. The proofs of the main results
are given in Section~\ref{secproofs}; Section~\ref{secunknowndensity}
discusses an unknown design density, and the proofs of lemmas are moved
to the supplementary material [\citet{suppA}].

\section{Main results}
\label{secmain-results}
In this section we motivate and explain our procedure and prove the
local and global oracle inequalities. Then we apply these results to
adaptive estimation over a collection of H\"older classes\vadjust{\goodbreak} (pointwise
estimation) and over a collection of Nikol'skii classes (estimating the
entire function with the accuracy of an estimator measured under
the~$L_r$~risk).

\subsection{Oracle approach}
\label{subsecoracle-approach}

Let $\mathcal{K}\dvtx \mathbb{R}\to\mathbb{R}$ be a function (kernel)
satisfying \mbox{$\int\mathcal{K}=1$}. With any such $\mathcal{K}$, any
$z\in\mathbb{R}$, $h\in(0,1]$ and any $f\in\mathbb{F}(\beta_0,M)$, we
associate the quantity
\[
\Delta_{\mathcal{K},f}(h,z)=\sup_{\delta\le h}\biggl\llvert
\frac
{1}{\delta}\int \mathcal{K} \biggl( \frac{u-z}{\delta} \biggr)
\bigl[f(u)-f(z) \bigr]\,\mathrm{d}u\biggr\rrvert.
\]
Note that the kernel smoother $ \delta^{-1} \int\mathcal{K}
([u-z]/\delta )f(u)\,\mathrm{d}u$ can be understood as an approximation
of the function $f$ at the point $z$. Thus,
$\Delta_{\mathcal{K},f}(h,z)$ is a monotonous approximation error
provided by this kernel smoother. In particular, under Assumption
\ref{asstechnical}, we have $ \Delta_{\mathcal{K},f}(h,z)\to0 $ as $
h\to0 $.

In what follows, $\|\mathcal{K}\|_p,  1\le p\le\infty$, denotes the
$L_p$ norm of $\mathcal{K}$ and we will assume that the kernel
$\mathcal{K}$ satisfies the following condition.

\begin{assumption}\label{assassumption-on-kernel}
(1)~$\operatorname{supp}(\mathcal{K})\subseteq[-1/2,1/2]$, $\int
\mathcal{K}=1$, $\mathcal{K}$ is symmetric;

(2)~there exists $Q>0$ such that $  |\mathcal{K}(u)-\mathcal {K}(v)
|\le Q|u-v|$ $\forall u,v\in\mathbb{R}$.
\end{assumption}

\subsubsection{Oracle estimator}
For any $y\in\mathbb{R}$, denote the Hardy--Littlewood maximal function
of $\Delta_{\mathcal{K},f}(h,\cdot)$ [see, e.g.,
\citet{WheedenZygmund1977}] by
\[
\Delta^*_{\mathcal{K},f}(h,y)=\sup_{a>0} \frac{1}{2a}
\int_{y-a}^{y+a}\Delta_{\mathcal{K},f}(h,z)\,\mathrm{d}z.
\]
Clearly, $  \Delta^{*}_{\mathcal{K},f}(h,\cdot)\ge\Delta_{\mathcal
{K},f}(h,\cdot) $ for any $ f\in\mathbb{F}(\beta_0,M) $.
Now, let us define the oracle estimator. For any $y\in\mathbb{R}$ and
$h_{\min}$ defined in~(\ref{eqbandwidts}), set
%
\begin{equation}
\label{deforaclebandwidth} h^*_{\mathcal{K},f}(y)= \sup \bigl\{h\in [h_{\min}, 1
] | \sqrt{nh} \Delta^{*}_{\mathcal{K},f}(h,y)\le \| \mathcal{K}
\|_{\infty} \sqrt{\ln(n)} \bigr\}.
\end{equation}
Note that $\Delta^{*}_{\mathcal{K},f}(h,\cdot)\le M\|\mathcal{K}\|
_1h^{\beta_0}$ for any $f\in\mathbb{F}(\beta_0,M)$ and any $h>0$.
Hence, $
\sqrt{nh_{\min}}\Delta^{*}_{\mathcal{K},f}(h_{\min},\cdot)\le\|
\mathcal{K}\| _1\sqrt{\ln(n)} $ for any $n\ge n_0$ in view of
(\ref{eqrestriction-on-sample-size}).
Next, Assumption \ref{assassumption-on-kernel}(2) implies that
$\Delta^{*}_{\mathcal{K},f}(\cdot,y)$ is continuous, hence,
%
\begin{eqnarray}
\label{eq1deforaclebandwidth} &\mbox{either}&\quad \sqrt{ nh^*_{\mathcal{K},f}(y)} \Delta
^{*}_{\mathcal{K},f} \bigl( h^*_{\mathcal{K},f}(y),y \bigr) = \|\mathcal{K}\|_{\infty}\sqrt{\ln(n)},
\\
\label{eq2deforaclebandwidth} &\mbox{or}&\quad\sqrt{nh} \Delta^{*}_{\mathcal{K},f}(h,y)
\le \| \mathcal{K}\|_{\infty}\sqrt{\ln(n)}\qquad\forall h\in[h_{\min},1].
\end{eqnarray}
Here we have also used that $\| \mathcal{K}\|_{1}\le\| \mathcal{K}\|
_{\infty}$ in view of Assumption~\ref{assassumption-on-kernel}(1).

The quantity similar to $h^*_{\mathcal{K},f}$ first appeared in
\citet{LepMamSpok97} for estimating univariate functions with
inhomogeneous smoothness. Some years later, this idea was further
developed for multivariate function estimation; see
\citet{KLP2001}, \citet{GoldLep2008} and the more detailed
discussion of the oracle approach therein. Following their lead, we
advance it for the estimation under the single-index constraint. The
basic idea behind our selection rule is simple.

For any $(\theta,h)\in\mathbb{S}^{1}\times[h_{\min}, 1]$, define the
matrix
\[
E_{(\theta,h)}=\pmatrix{ h^{-1}\theta_1 &
h^{-1}\theta_2
\cr
-\theta_2 &
\theta_1},\qquad\det (E_{(\theta, h)} )=h^{-1}
\]
and consider the family of kernel estimators with $ K(u,v)=\mathcal
{K}(u)\mathcal{K}(v) $ so that
\[
\label{familyofestimators} \mathcal F = \Biggl\{ \widehat F_{(\theta, h)} (\cdot) =
\frac{\det (E_{(\theta, h)} )}{n}\sum_{i=1}^n
\frac{K
(E_{(\theta, h)}(X_i- \cdot) )}{g(X_i)} Y_i, (\theta,h)\in \mathbb{S}^{1}
\times[h_{\min}, 1] \Biggr\}.
\]
We remark that Assumptions \ref{assassumption-on-design} and
\ref{assassumption-on-kernel}(1) assure well-definiteness of $\widehat
F_{(\theta, h)}$ because $ K (E_{(\theta, h)}(x-t) )=0$  $\forall
x\in[-3/2,3/2]^2 $ and $ \forall t\in[-1/2,1/2]^2 $.

The\vspace*{-2pt} choice $\theta=\theta^*$ and $h=h^*:=h^*_{\mathcal
{K},f}(t^{\top}\theta^*)$ leads to the \textit{oracle estimator} $\widehat
F_{(\theta^*, h^*)}$. Note that $\widehat F_{(\theta^*, h^*)}$ is not
an estimator in the usual sense because it depends on the function $F$
to be estimated [more precisely, on $(f,\theta^*)$ which
determines~$F$]. The meaning of $ \widehat F_{(\theta^*, h^*)} $ is
explained by the following result based
on the straightforward application of Rozenthal's inequality.

\begin{proposition}
\label{proprisk-of-oracle-estimator} For any $
(f,\theta^*)\in\mathbb{F}(\beta_0,M)\times\mathbb{S}^{1}$,  $r \ge1 $
and $n\ge n_0$,
\[
\mathcal R_{r,t}^{(n)} (\widehat{F}_{(\theta^*, h^*)}, F ) \le
c \ln^{1/2}(n) \bigl[nh^*_{\mathcal{K},f}\bigl(t^{\top}\theta^*
\bigr) \bigr]^{-1/2}\qquad \forall t\in[-1/2,1/2]^{2},
\]
where $c>0$ is a numerical constant independent of $n$.
\end{proposition}

This result indicates that the ``oracle'' knows the exact value of
$\theta^*$ and the optimal, up to $\ln(n)$, trade-off $h^*$ between the
approximation error induced by $\Delta^{*}_{\mathcal
{K},f}(h^{*},\cdot)$ and the stochastic error of the kernel estimator
from $\mathcal{F}$ with bandwidth $h^*$. It explains why the ``oracle''
chooses the ``estimator'' $\widehat{F}_{(\theta^*, h^*)}$. Below we
propose a ``real,'' based on the observation, estimator $\widehat{F}$,
which mimics the oracle---for any
$(f,\theta^*)\in\mathbb{F}(\beta_0,M)\times\mathbb{S}^{1}$, $r\ge1$ and
$n\ge n_0$,
\[
\mathcal R_{r,t}^{(n)} (\widehat{F}, F ) \le
c^\prime \ln^{1/2}(n) \bigl[ nh^*_{\mathcal{K},f}
\bigl(t^{\top}\theta^*\bigr) \bigr]^{-1/2}\qquad\forall t
\in[-1/2,1/2]^{2},
\]
where $c^\prime$ is an absolute constant independent of~$n$ and the
underlying function $F$. The latter result is a local oracle
inequality. The construction of the estimator $\widehat{F}$ is based on
the data-driven selection from the family $\mathcal{F}$.

\subsubsection{Selection rule}
\label{secselection-rule}

For any $\theta,\nu\in\mathbb{S}^{1}$ and any $h\in[h_{\min},1]$,
define
\[
\widebar{E}_{(\theta,h)(\nu,h)}=\pmatrix{\displaystyle\frac{(\theta_1+\nu_1)}{2h(1+|\nu^\top\theta|)} &\displaystyle\frac{(\theta_2+\nu
_2)}{2h(1+|\nu^\top\theta|)}
\vspace*{6pt}\cr
\displaystyle -\frac{(\theta_2+\nu_2)}{2(1+|\nu^\top\theta|)} & \displaystyle\frac{(\theta_1+\nu _1)}{2(1+|\nu^\top\theta|)}},
\]
where
\[
E_{(\theta, h)(\nu,h)}=\cases{ \widebar{E}_{(\theta, h)(\nu,h)},&\quad$\nu^\top
\theta\ge0$,
\vspace*{3pt}\cr
\widebar{E}_{(-\theta, h)(\nu,h)},&\quad$\nu^\top\theta<
0$,} \qquad\frac{1}{4h} \le\det (E_{(\theta,h)(\nu,h)} )\le\frac{1}{2h}.
\]
%
A kernel estimator associated with the matrix $E_{(\theta, h)(\nu, h)}$
is defined by
%
\begin{equation}
\label{pseudo-estimator} \widehat F_{(\theta, h)(\nu, h)}(\cdot) = \frac{ \det (E_{(\theta, h)(\nu,h)} )}{n} \sum
_{i=1}^n \frac{K( E_{(\theta, h)(\nu,h)}(X_i- \cdot)
)}{g(X_i)}
Y_i.
\end{equation}
The definition of $\widehat F_{(\theta,
h)(\nu, h)}$ is legitimate because $ K (E_{(\theta,
h)}(x-t) )=0 $ $ \forall x\in[-5/2,5/2]^2 $ and $ \forall
t\in[-1/2,1/2]^2 $.

For any $u_1,u_2\in\mathbb{R}$, set $K_{\mathfrak
{h}}(u_1,u_2)=\mathfrak{h}^{-2}\mathcal{K}(u_1/\mathfrak{h})\mathcal
{K}(u_2/\mathfrak{h})$ and define
\[
\widehat{F}(v)=n^{-1}\sum_{i=1}^n
g^{-1}(X_i)K_\mathfrak{h}(X_i-
v)Y_i,\qquad \widehat{F}_\infty=2\|\widehat{F}
\|_\infty+2C_5(n),
\]
where $ \|\widehat{F}\|_\infty=\sup_{v\in[-5/2,5/2]^2 } |\widehat
{F}(v) | $ and $\mathfrak{h}$ is defined in (\ref{eqbandwidts}). Put
also
\[
\operatorname{TH}(\eta)=2 \bigl[\|\mathcal{K}\|^{2}_{\infty}
\sqrt {\ln(n)}+\widehat{F}_\infty C_1(n)+C_2(n)
\bigr](\eta n)^{-1/2},\qquad \eta\in(0,1].
\]
The quantities $C_1(n)$, $C_2(n)$ and $C_5(n)$ are listed in Section~\ref{secconstants}.

Set $\mathcal H_n= \{h_k=2^{-k},  k \in \mathbb{N}^0 \}\cap
[2^{-1}h_{\min},1 ]$ and let for any $\theta\in\mathbb{S}^{1}$ and $h
\in\mathcal H_n $,
\begin{eqnarray*}
R^{(1)}_{t}(\theta, h) &=& \sup_{\eta\in\mathcal H_n\dvtx   \eta\le h}
\Bigl[ \sup_{\nu\in\mathbb{S}^{1}} \bigl|\widehat F_{(\theta, \eta)(\nu, \eta)}(t) - \widehat
F_{(\nu, \eta
)}(t) \bigr| -\operatorname{TH}(\eta) \Bigr]_+,
\\
R^{(2)}_{t}(h)&=&\sup_{\eta\in\mathcal H_n\dvtx   \eta\le h} \Bigl[\sup
_{\theta\in\mathbb{S}^{1}} \bigl|\widehat F_{(\theta,h)}(t)- \widehat
F_{(\theta,
\eta
)}(t) \bigr|-\operatorname{TH}(\eta) \Bigr]_+.
\end{eqnarray*}
Subsequently, define $(\hat{\theta},\hat{h})$ as a solution of the
following minimization problem:
%
\begin{eqnarray}
\label{eqselection-rule} &&R^{(1)}_{t}(\hat{\theta},
\hat{h})+R^{(2)}_{t}(\hat{h})+\operatorname{TH}(
\hat{h})
\nonumber\\[-8pt]\\[-8pt]
&&\qquad =\inf_{(\theta,h)\in\mathbb{S}^{1}\times\mathcal H_n} \bigl[R^{(1)}_{t}(\theta,
h)+R^{(2)}_{t}( h)+\operatorname{TH}(h) \bigr].\nonumber
\end{eqnarray}
Then our final estimator is $ \widehat F (t)=\widehat F_{(\hat \theta,
\hat h)}(t)$, where $ (\hat\theta, \hat h) $ is obtained by minimizing
(\ref{eqselection-rule}).

\begin{remark}
\label{remmeasurability} We note that
Assumption~\ref{assassumption-on-kernel}(2) guarantees that all random
fields involved in the description of selection rule
(\ref{eqselection-rule}) are continuous on~$\mathbb{S}^{1}$. Moreover,
the set $\mathcal H_n$ is finite. Thus, $(\hat{\theta},\hat{h})$ is $
\{ (X_i,Y_i)\}_{i=1}^n$-measurable and $(\hat{\theta},\hat{h})\in
\mathbb{S}^{1}\times\mathcal H_n$ [see \citet{jennrich}].
\end{remark}

\begin{remark}
\label{remdifficulties} Our selection rule (\ref{eqselection-rule}) is
defined in the case $d=2$. The main difficulty in extending it to $d>2$
consists in the construction of the matrix $E_{(\theta,h)(\nu,h)}$ for
any vectors $\theta,\nu\in\mathbb{S}^{d-1}$. Indeed, analyzing the
proof of Theorem \ref{thlocal-oracle-inequality}, we remark that the
following properties should be fulfilled:
\[
E_{(\theta,h)(\nu,h)}\in\mathcal E_{a,A}, \qquad E_{(\theta, h)(\nu,h)}=\pm
E_{(\nu, h)(\theta,h)}\qquad \forall\theta,\nu\in\mathbb {S}^{d-1}, \forall h\in
\mathcal H_n,
\]
where the class of matrices $\mathcal E_{a,A}$ is defined in
(\ref{eqdefmatrixclass}). If $d=2$, these requirements hold. However,
we were not able to construct a class of matrices obeying latter
restrictions in the dimension strictly larger than 2. Note,
nevertheless, that if such a class would be found, our results could be
extended to $d>2$ without any additional consideration.
\end{remark}

\subsubsection{Local and global oracle inequalities} We reinforce
restriction~(\ref{eqrestriction-on-sample-size}) on the minimal sample
size $n$. Let $n_1\ge1$ be defined as follows:
%
\begin{equation}
\label{eqrestriction2-on-sample-size} n_1=\inf \bigl\{m\in\mathbb{N}^*\dvtx  \bigl(n
\mathfrak{h}^2 \bigr)^{-1/2}C_3(n)\le1/2\ \forall n\ge m \bigr\},
\end{equation}
where $\mathfrak{h}$ is defined in (\ref{eqbandwidts}) and $C_3(n)$ is
given at the beginning of Section~\ref{secconstants}. All our results
below will be proved under the condition $n\ge n_0\vee n_1$.

First, we note that $n_1$ is well-defined since 
$ (n\mathfrak{h}^2 )^{- 1/2}C_3(n)\to0$ as $n\to\infty$. Next, contrary
to restriction (\ref{eqrestriction-on-sample-size}) that relates the
sample size $n$ to the quantities $ \beta_0 $ and $ M $ from Assumption
\ref{asstechnical}, restriction (\ref{eqrestriction2-on-sample-size})
links the minimal value of $n$ with the quantity $\underline{g}$
appearing in Assumption \ref{assassumption-on-design}.

\begin{theorem}
\label{thlocal-oracle-inequality} For any
$(f,\theta^*)\in\mathbb{F}(\beta_0, M)\times\mathbb{S}^{1}$, $r \ge1$
and $n\ge n_0\vee n_1$,
\[
\mathcal R_{r,t}^{(n)} (\widehat{F}_{(\hat{\theta},
\hat{h})}, F ) \le
c_{1} \biggl[\frac{ \ln(n)}{nh^*_{\mathcal{K},f}(t^{\top}\theta
^*)} \biggr]^{1/2} +
c_{2}n^{-1/2}\qquad\forall t\in[-1/2,1/2]^{2}.
\]
\end{theorem}

The constants $c_{1}$ and $c_{2}$ are independent of $n$ and $F$ and
their explicit expressions can be extracted from the proof of the
theorem.

As already mentioned, the global oracle inequality is obtained by
integrating the local oracle inequality. Indeed, for any $ r \ge1 $,
using Jensen's inequality and Fubini's theorem, we have $ \mathcal
R_r^{(n)} (\widehat F, F) \le \|\mathcal R_{r,\cdot}^{(n)} (\widehat F,
F) \|_r $ so
\[
\mathcal R_r^{(n)} (\widehat F, F) \le c_{1}
\biggl\{ \int_{[-1/2, 1/2]^{2}} \biggl[\frac{
\ln(n)}{nh^*_{\mathcal{K},f}(t^{\top}\theta^*)}
\biggr]^{r/2}\,\mathrm{d}t \biggr\}^{1/r} + c_{2}n^{-1/2}.
\]
Integration by substitution yields
\[
\int_{[-1/2, 1/2]^{2}} \biggl[\frac{ \ln(n)}{nh^*_{\mathcal{K},f}(t^{\top}\theta^*)} \biggr]^{r/2}\,\mathrm{d}t \le \int_{-1/2}^{1/2} \biggl[
\frac{ \ln(n)}{nh^*_{\mathcal{K},f}(z)} \biggr]^{r/2}\,\mathrm{d}z,
\]
that leads to the following bound.

\begin{theorem}
\label{thglobal-oracle-inequality} For any
$(f,\theta^*)\in\mathbb{F}(\beta_0,M)\times\mathbb{S}^{1}$, $r\ge1$ and
$n\ge n_0\vee n_1$,
\[
\mathcal R_r^{(n)} (\widehat{F}_{(\hat{\theta}, \hat{h})}, F )\le
c_{1} \biggl\|\frac{
\ln(n)}{nh^*_{\mathcal{K},f}} \biggr\|^{1/2}_{r/2}+c_{2}n^{-1/2}.
\]
\end{theorem}

\subsubsection{\texorpdfstring{Extension to the case $\theta^*\notin\mathbb{S}^{1}$}
{Extension to the case theta* not in S1}}
\label{secextensions}

Define $f_{\theta^*}(t)=f (|\theta^*|_2t )$, $\vartheta^*=\break \theta ^*
/|\theta^*|_2$ and let $F_{\theta^*}(t):=f_{\theta^*} (t^\top
\vartheta^* )$. Obviously, for all $t\in \mathbb{R}^2$ we have\break
$f_{\theta^*} (t^\top\vartheta^* )=f
(t^\top\theta ^* )$ that implies $F_{\theta^*}(
\cdot)\equiv F(\cdot)$
so the estimation of $F$ is equivalent to the estimation of $F_{\theta
^*}$. Because $\vartheta\in\mathbb{S}^{1}$,
Theorems~\ref{thlocal-oracle-inequality} and
\ref{thglobal-oracle-inequality} are applicable. To this end, it
suffices to replace $f$ by $f_{\theta^*}$ in the definition of
$h^*_{\mathcal{K},f}(\cdot)$. In general, however, there is no
universal\vspace*{1pt} way of expressing $h^*_{\mathcal{K},f_{\theta^*}}(\cdot)$ via
$h^*_{\mathcal {K},f}(\cdot)$, although in particular cases, mainly in
adaptive estimation over classes of smooth functions, it is often
possible.

\subsection{Adaptive estimation}\label{subsecadaptive-estimation}
In this section we first apply the local oracle inequality given in
Theorem \ref{thlocal-oracle-inequality} to the problem of pointwise
adaptive estimation over a collection of H\" older classes. Next, we
study adaptive estimation under the $L_r$ losses over a collection of
Nikol'skii classes. The corresponding result is deduced from the global
oracle inequality proved in Theorem~\ref{thglobal-oracle-inequality}.

Assume throughout this section that the kernel $\mathcal{K}$ obeys
additionally Assumption \ref{ass2assumption-on-kernel} below; we then
introduce the following notation: for any $a>0$, let $m_a $ be the
maximal integer strictly less than $a$.

%
\begin{assumption}
\label{ass2assumption-on-kernel} There exists $\mathbf{b}>0$ such that
\[
\int z^{j}\mathcal{K}(z)\,\mathrm{d}z=0\qquad \forall j=1,\ldots,m_{\mathbf{b}}.
\]
\end{assumption}

\subsubsection{Pointwise adaptive estimation}
We start with some definitions.

%
\begin{definition}\label{defholder-class}
Let $ \beta>0 $ and $ L>0 $. A function $\ell\dvtx
\mathbb{R}\to\mathbb{R}$ belongs to the H\"older class
$\mathbb{H}(\beta,L)$ if $\ell$ is $m_\beta$-times continuously
differentiable, $\|\ell^{(m)}\|_\infty\le L$ for all $m\le m_\beta$,
and
\[
\bigl\llvert \ell^{(m_\beta)}(u+h)-\ell^{(m_\beta)}(u)\bigr\rrvert \le L
h^{\beta-m_\beta}\qquad\forall u\in\mathbb{R}, h>0.
\]
\end{definition}

The aim is to estimate the function $F(t)$ at a given point $t\in
[-1/2,1/2]^2$ under the assumption that $F\in\mathbb{F}(\mathbf{b}
):=\bigcup_{\beta\le\mathbf{b}}\bigcup_{L>0}\mathbb{F}_2(\beta,L)$,
where
\[
\mathbb{F}_d(\beta,L)= \bigl\{F\dvtx \mathbb{R}^d\to
\mathbb{R} | F(z)=f \bigl(z^{\top
}\theta \bigr), f\in\mathbb{H}(\beta,L),
\theta\in\mathbb {S}^{d-1} \bigr\},\vadjust{\goodbreak}
\]
the constant $\mathbf{b}$ is from Assumption
\ref{ass2assumption-on-kernel}, and $d\ge2$ is the dimension. We will
see that $\mathbf{b}$ can be an arbitrary number but it must be chosen
{a priori}.\vspace*{-1pt}

%
\begin{theorem}
\label{thpointwise-adaptation} Let $\mathbf{b}>0$ be fixed; and let
additionally Assumptions \ref{assassumption-on-kernel} and
\ref{ass2assumption-on-kernel} hold. Then, for any
$\beta\le\mathbf{b}$, $L>0$,  $r \ge1 $ and $t\in[-1/2,1/2]^{2}$,
\[
\sup_{F\in\mathbb{F}_2(\beta,L)}\mathcal R_{r,t}^{(n)} (\widehat{F}_{(\hat{\theta}, \hat{h})}, F ) \le \varkappa_1\psi_n(
\beta,L),
\]
where $\psi_n(\beta,L)=L^{1/(2\beta+1)} [n^{-1}\ln (n)
]^{\beta/(2\beta+1)}$ and $\varkappa_1$ is independent of $n$.\vspace*{-1pt}
\end{theorem}

The proof of the theorem is based on the evaluation of the uniform over
$\mathbb{H}_d(\beta,L)$ lower bound for $ h^*_{\mathcal{K},f}(\cdot) $
and on the application of Theorem~\ref{thlocal-oracle-inequality}. We
note that a similar upper bound for the minimax risk appeared in
\citet{GoldLep2008} in the framework of Gaussian white noise
model, but the estimation procedure used there is different from our
selection rule.

The main question, however, is if $\psi_n(\beta,L)$ coincides with the
minimax rate for any given value of $\beta$ and $L$? To answer it, we
need some additional assumptions on the densities of the noise variable
$\varepsilon_1$ and design variable $X_1$.\vspace*{-1pt}

\begin{assumption}
\label{ass2on-noise} There exist $\mathfrak{q},\mathfrak{Q}>0$ such
that, for any $\upsilon_1,\upsilon_2\in[-\mathfrak{q},\mathfrak{q}]$,
\[
\int_{\mathbb{R}} p(y+\upsilon_1)p(y+
\upsilon_2)p^{-1}(y) \,\mathrm{d}y\le 1+\mathfrak{Q} |
\upsilon_{1}\upsilon_{2} |.\vspace*{-1pt}
\]
\end{assumption}

It is easy to see that the density of the normal law
$\mathcal{N}(0,\sigma^2)$,  $\sigma^2>0 $, obeys the aforementioned
assumption.
In general, this assumption is fulfilled if the density $p$~is regular
and decreases rapidly at infinity. More precisely, if the Fisher
information corresponding to the density $p$ is finite and the
function $\int [p^\prime(y+\cdot) ]^2p^{-1}(y) \,\mathrm{d}y$ is
continuous at zero, Assumption \ref{ass2on-noise} is verified.\vspace*{-1pt}

\begin{assumption}
\label{ass2on-design} There exist $\mathfrak{g}>0$ and $\varpi> 1$ such
that, for all $ x\in\mathbb{R}^d$,  $g(x)\le (1+|x|_2^{\varpi}
)^{-1}\mathfrak{g}$.  Here $|\cdot|_2$ is the Euclidean vector norm on
$\mathbb{R}^d$.\vspace*{-1pt}
\end{assumption}

We remark that the imposed assumption is very weak and holds for the
majority of probability distributions used in statistical applications.\vspace*{-1pt}

\begin{theorem}
\label{thpointwise-adaptation-lower} Let Assumptions \ref{ass2on-noise}
and \ref{ass2on-design} be fulfilled. Then, for any $ t
\in[-1/2,1/2]^d$,  $d\ge2 $,  $r \ge1 $,  $\beta, L>0 $, and any
$n\in\mathbb{N}^*$ large enough,
\[
\inf_{\widetilde{F}}\sup_{F\in\mathbb{F}_d(\beta,L)}\mathcal
R_{r,t}^{(n)} (\widetilde{F}, F ) \ge\varkappa_2
\psi_n(\beta,L),
\]
where the infimum is over all possible estimators. Here $\varkappa_2$
is a numerical constant independent of $n$ and $L$, and $
\psi_n(\beta,L) $ is defined in Theorem \ref{thpointwise-adaptation}.\vspace*{-1pt}
\end{theorem}

To the best of our knowledge, this lower bound is new. It is worth
mentioning that Assumption \ref{ass2on-noise} is close to being
necessary. One can give examples where this condition does not hold and
Theorem \ref{thpointwise-adaptation-lower} is not true anymore.\vadjust{\goodbreak}

Theorems \ref{thpointwise-adaptation} and
\ref{thpointwise-adaptation-lower} indicate that the estimator
$\widehat {F}_{(\hat{\theta}, \hat{h})}$ is minimax adaptive with
respect to the collection $ \{\mathbb{F}_d(\beta,L), \beta\le\mathbf
{b}, L>0 \}$. As already mentioned, this result is quite surprising.
Indeed, if, for example, $\theta=(1,0)^{\top}$, that is, is known, then
$\mathbb{F}(\beta,L)=\mathbb{H}(\beta,L)$, and the considered
estimation problem reduces to estimation of $f$ at a point in the
univariate regression model. As it is shown in~\citet{Gaiffas}, an
adaptive estimator over $ \{ \mathbb{H}(\beta,L), \beta\le\mathbf{b},
L>0 \}$ does not exist and a price for adaption appears. The latter
means that the asymptotic bound on the minimax risk provided by the
adaptive estimator differs from the minimax rate of convergence by some
factor. This factor for the majority of known results is $\ln(n)$.

In addition, we would like to note that the assertion of
Theorem~\ref{thpointwise-adaptation-lower} is proved for arbitrary
dimension.

\subsubsection{Adaptive estimation under the $ L_r $ losses}
We begin by defining the relevant functional classes.

\begin{definition}\label{defnikolskii-class}
Let $ \beta>0
$,  $L>0 $ and $p\in[1,\infty)$ be fixed. A function $\ell\dvtx
\mathbb{R}\to\mathbb{R}$ belongs to the Nikol'skii class
$\mathbb{N}_p(\beta,L)$ if $\ell$ is $m_\beta$-times continuously
differentiable and
\begin{eqnarray*}
\biggl( \int_{\mathbb{R}}\bigl\llvert \ell^{(m)}(t)
\bigr\rrvert ^p \,\mathrm{d}t \biggr)^{1/p}&\le& L\qquad \forall m=0,\ldots, m_\beta,
\\
\biggl( \int_{\mathbb{R}}\bigl\llvert \ell^{(m_\beta)}(t+h)-
\ell ^{(m_\beta
)}(t)\bigr\rrvert ^p \,\mathrm{d}t
\biggr)^{1/p} &\le& Lh^{\beta-m_\beta}\qquad\forall h>0.
\end{eqnarray*}
\end{definition}

It is also assumed that $\mathbb{N}_p(\beta,L)=\mathbb {H}(\beta,L)$ if
$p=\infty$.

Here, the target of estimation is the entire function $F$ under the
assumption that $ F
\in\mathbb{F}_p(\mathbf{b}):=\bigcup_{\beta\le\mathbf{b}}\bigcup_{L>0}\mathbb{F}_{2,p}(\beta,L)$,
where
\[
\mathbb{F}_{d,p}(\beta,L)= \bigl\{F\dvtx
\mathbb{R}^d\to\mathbb{R} | F(z)=f \bigl(z^{\top}\theta
\bigr), f\in\mathbb{N}_p(\beta,L), \theta\in \mathbb{S}^{d-1}
\bigr\}.
\]
Let us briefly discuss the applicability of Theorem
\ref{thglobal-oracle-inequality} requiring $f\in\mathbb{F}(\beta_0,M)$.
To this end, we assume that $\beta p> 1$. The latter assumption is
standard for estimating functions with inhomogeneous smoothness [see,
e.g., \citet{DJKP}, \citet{LepMamSpok97},
\citet{KLP2008}]. If $\beta p> 1$, the embedding
$\mathbb{N}_p(\beta,L)\subset\mathbb{H}(\beta-1/p,cL)$ with an absolute
constant $c>0$ guarantees that $f\in\mathbb{F}(\beta_0,M)$ with $
\beta_0=\beta-1/p $ and $ M=cL $.

\begin{theorem}\label{thglobal-adaptation}
Let $\mathbf{b}>0$ be fixed, and let
Assumptions \ref{assassumption-on-kernel} and
\ref{ass2assumption-on-kernel} hold. Then, for any $L>0$, $p>1$,
$p^{-1}<\beta\le\mathbf{b}$ and $r\ge1$,
\[
\sup_{F\in\mathbb{F}_{2,p}(\beta,L)}\mathcal R_{r}^{(n)} (\widehat
{F}_{(\hat{\theta}, \hat{h})}, F ) \le \varkappa_3\varphi_n(
\beta,L,p),
\]
where $\varkappa_3$ is independent of $n$, and
\[
\varphi_n(\beta,L,p)=\cases{ L^{1/(2\beta+1)} \bigl(n^{-1}
\ln(n) \bigr)^{\beta/(2\beta+1)},
\vspace*{2pt}\cr
\qquad (2\beta+1)p>r,
\vspace*{2pt}\cr
L^{1/(2\beta+1)}\bigl(n^{-1}\ln(n) \bigr)^{\beta/(2\beta+1)}\ln^{1/r}(n),
\vspace*{2pt}\cr
\qquad (2\beta+1)p=r,
\vspace*{2pt}\cr
L^{(1/2-1/r)/(\beta-1/p+1/2)} \bigl(n^{-1}\ln(n)
\bigr)^{(\beta-1/p+1/r)/(2\beta-2/p+1)},
\cr
\qquad (2\beta+1)p< r.}
\]
\end{theorem}

Note that $\mathbb{F}_{2,p}(\beta,L)\supset\mathbb{N}_{p}(\beta,L)$.
Indeed, the class $\mathbb{N}_{p}(\beta,L)$ can be viewed as a class of
functions $F$ satisfying $F(\cdot)=f(\theta^\top\cdot)$ with
$\theta=(1,0)^{\top}$. Then, the problem of estimating such (2-variate)
functions reduces to the estimation of univariate regression functions.

There are at least two observations arising in view of the latter
remark. First, the upper bound of Theorem \ref{thglobal-adaptation}
generalizes the results for the univariate regression
[\citet{DJKP}, \citet{DelyonJuditski1996},
\citet{Baraud}, \citet{KP}, \citet{Kulik},
\citet{Zhang}] in several directions. In particular, the majority
of the papers treat the Gaussian errors or the errors having
exponential moment. An exception is \citet{Baraud}, where some
results are obtained under a very weak assumption on the noise (weaker
than our Assumption~\ref{assassumption-on-noise}). Nevertheless, these
results are available only if $p=r=2$.

Next, the rate of convergence for the latter problem, which can be
found in \citet{Chesneau}, is also the lower bound for the minimax
risk defined on $\mathbb{F}_{2,p}(\beta,L)$. With the proviso that
$\beta p>1$, the rate of convergence is given by
\[
\phi_n(\beta,L,p)=\cases{ L^{1/(2\beta+1)} n^{-\beta/(2\beta+1)},
\cr
\qquad (2\beta+1)p> r,
\vspace*{2pt}\cr
L^{1/(2\beta+1)} \bigl(n^{-1}\ln(n)
\bigr)^{\beta/(2\beta +1)},
\vspace*{2pt}\cr
\qquad (2\beta+1)p=r,
\vspace*{2pt}\cr
L^{(1/2-1/r)/(\beta-1/p+1/2)}
\bigl(n^{-1}\ln(n) \bigr)^{(\beta-1/p+1/r)/(2\beta
-2/p+1)},
\cr
\qquad (2\beta+1)p< r.}
\]
The minimax rate of convergence in the case $(2\beta+1)p=r$ is not
known, hence, the rate presented in the middle line above is only the
lower asymptotic bound for the minimax risk.

Thus, the proposed estimator $\widehat{F}_{(\hat{\theta}, \hat{h})}$ is
adaptive whenever $(2\beta+1)p<r$. In the case $(2\beta+1)p\ge r$, we
loose only a logarithmic factor with respect to the optimal rate and,
as mentioned in the \hyperref[SectionIntro]{Introduction}, the construction of an adaptive
estimator over the collection $ \{\mathbb{F}_{2,p}(\beta,L),  \beta>0,
L>0 \}$ in this case remains an open problem. In view of the latter
remark, we conjecture that the presented lower bound is correct and,
therefore, the upper bound result has to be improved.

\section{Proofs}\label{secproofs}
We now list the quantities that are involved in the description of the
selection rule that led to the adaptive estimator $\widehat
F_{(\hat\theta, \hat h)}$.

\subsection{Important quantities}
\label{secconstants} Let $ \tau= (\Omega^{-1}(4r+1)\ln(n) )^{1/\omega}
$. Set
\begin{eqnarray*}
c_1(n) &=& 730\ln{ \bigl(16n^2
\underline{g}^{-1/2} [12Q+\sqrt {2} ] \bigr)}+ 8r\ln(n)+394,
\\
c_2(n) &=& 730\ln{ \bigl(16n^2\tau\underline{g}^{-1/2}
[12Q +\sqrt{2} ] \bigr)}+ 8r\ln(n)+394,
\\
c_3(n)&=&365\ln{ \bigl(5n^{2}Q\underline{g}^{-1/2}
\bigr)}+8r\ln(n)+197,
\\
c_4(n)&=&365\ln{ \bigl(5n^2\tau
Q \underline{g}^{-1/2} \bigr)}+8r\ln(n)+197.
\end{eqnarray*}
With $ \mathfrak{h} $ given in (\ref{eqbandwidts}) and $
\sigma^{2}=\sup_{p\in\mathfrak{P}}\int_{\mathbb{R}}x^2p(x)\,\mathrm{d}x
$, we define
\begin{eqnarray*}
C_1(n)&=&2\sqrt{2}\underline{g}^{-1/2}\|\mathcal{K}\|
^2_\infty\sqrt {c_1(n)}+ (8/3)c_1(n)
\bigl(\ln(n) \bigr)^{-(2+\omega)/(2\omega)}\underline {g}^{-1}\|\mathcal{K}
\|^2_\infty,
\\
C_2(n)&=&2\sqrt{2}(\sigma\vee1)\underline{g}^{-1/2}\|
\mathcal{K}\| ^2_\infty\sqrt{c_2(n)}
\\
&&{}+ (8/3) c_2(n) \bigl(\ln(n) \bigr)^{-1/2}
\underline{g}^{-1}\| \mathcal{K} \|^2_\infty \bigl(
\Omega^{-1}(4r+1) \bigr)^{1/\omega},
\\
C_3(n)&=&2\sqrt{2}\underline{g}^{-1/2}\|\mathcal{K}\|
^2_\infty\sqrt {c_3(n)}+ (8/3)
\underline{g}^{-1}\|\mathcal{K}\|^2_\infty
c_3(n) \bigl(n\mathfrak {h}^2 \bigr)^{-1/2},
\\
C_4(n)&=&2\sqrt{2}(\sigma\vee1)\underline{g}^{-1/2}\|
\mathcal{K}\| ^2_\infty\sqrt{c_4(n)}+ (8/3)\tau
c_4(n) \bigl(n\mathfrak{h}^2 \bigr)^{-1/2}
\underline {g}^{-1}\|\mathcal{K}\|^2_\infty,
\\
C_5(n)&=&\|\mathcal{K}\|^2_1+ \bigl(n
\mathfrak{h}^2 \bigr)^{-1/2}C_4(n) + 1/2.
\end{eqnarray*}
%
In spite of the cumbersome expressions, it is easy to see that
%
\begin{equation}
\label{eqasympt-bound-for-quantities} \qquad\sup_{n\ge3}\frac{C_i(n)}{\sqrt{\ln(n)}}=:C_i<
\infty,\qquad i=1,2,\qquad\sup_{n\ge3}C_5(n)=:C_5<
\infty.
\end{equation}

\subsection{\texorpdfstring{Proof of Theorem \protect\ref{thlocal-oracle-inequality}}
{Proof of Theorem 1}}
To begin with we present upper bounds for the approximation errors
of the estimators involved (Lemma \ref{lemboundsforbias}) and their
stochastic errors (Lemma \ref{lemgauss-on-matrices}). Lemma
\ref{leminegality-for-sup-norm} allows us to proceed without
knowledge of  $M$ from Assumption \ref{asstechnical}. The proofs
of the later two results are essentially based on Proposition 1 of
\citet{Lep2013}. The detailed proofs of these technical results are
moved to the supplementary material [\citet{suppA}].

\subsubsection{Auxiliary results}
\label{secaux-results-thlocal-oracle-inequality} For any
$\theta,\nu\in\mathbb{S}^{1}$ and $h\in [2^{-1}h_{\min},1 ]$, denote
\begin{eqnarray*}
S_{(\theta, h)(\nu, h)}(t) &=& \det (E_{(\theta, h)(\nu,h)} ) \int K
\bigl(E_{(\theta, h)(\nu,h)}(x-t)\bigr)F(x)\,\mathrm{d}x,
\\
S_{(\theta, h)}(t) &=& \det (E_{(\theta, h)} ) \int K\bigl(E_{(\theta, h)}(x-t)
\bigr)F(x)\,\mathrm{d}x.
\end{eqnarray*}
For ease of notation, we write $h^*_f=h^*_{\mathcal{K},f}(t^{\top
}\theta^*)$.

%
\begin{lemma}
\label{lemboundsforbias} Grant Assumption
\ref{assassumption-on-kernel}. Then, for any $\nu\in\mathbb{S}^{1}$ and
any bandwidths $\eta,h\in [2^{-1}h_{\min},1 ]$ satisfying $\eta\le
h\le2^{-1}h^*_{f}$, one has
\begin{eqnarray*}
\bigl\llvert S_{(\theta^*, h)(\nu, h)}(t)-S_{(\nu, h)}(t)\bigr\rrvert &\le& 2
\bigl(h^*_f\bigr)^{-1/2}\| \mathcal{K}\|^{2}_{\infty}
\sqrt{n^{-1}\ln(n)},
\\
\bigl\llvert S_{(\nu, h)}(t)-S_{(\nu, \eta)}(t)\bigr\rrvert &\le& 2
\bigl(h^*_f\bigr)^{-1/2}\| \mathcal{K}\|^{2}_{\infty}
\sqrt{n^{-1}\ln(n)},
\\
\bigl\llvert S_{(\theta^*, h)}(t) - F(t)\bigr\rrvert &\le&\bigl(h^*_f
\bigr)^{-1/2}\| \mathcal{K}\| _{\infty}\sqrt{n^{-1}
\ln(n)}.
\end{eqnarray*}
\end{lemma}

Let $ \mathcal E_{a,A} $ with $a \in(0,1], A\ge1$, be a set of
$2\times2$ matrices satisfying
%
%
\begin{equation}
\label{eqdefmatrixclass} \bigl|\det(E) \bigr|\le A,\qquad|E|_\infty\le(\sqrt{2a})^{-1}
\bigl|\det(E) \bigr|.
\end{equation}
Here $|E|_\infty= \max_{i,j} |E_{i,j}|$ denotes the matrix $\operatorname{sup}$
norm. Set, $\forall E\in\mathcal E_{a,A}$,
\begin{eqnarray*}
J(x,E)&=&\sqrt{ \bigl|\det(E) \bigr|}K \bigl(E(x-t) \bigr)g^{-1}(x),\qquad x\in
\mathbb{R}^2 
\end{eqnarray*}
and consider the following random fields defined on $\mathcal E_{a,A}$:
\begin{eqnarray*}
\eta_{n,t}(E) &=& n^{-1/2}\sum_{i=1}^n
\bigl\{ J(X_i,E)F(X_i)-\mathbb{E}^{(n)}_{X}
\bigl[J(X_i,E)F(X_i) \bigr] \bigr\},
\\
\xi_{n,t}(E) &=& n^{-1/2}\sum_{i=1}^n
J(X_i,E) \varepsilon_i.
\end{eqnarray*}
Denote finally by $\mathcal E_*$ the set of matrices $\mathcal E_{a,A}$
with $a=1/8$ and $A=h_{\min}^{-1}$. In what follows, we denote by $ \|
F \|_\infty= \sup_{ x\in[-5/2, 5/2]^2 } |F(x)| $.

\begin{lemma}
\label{lemgauss-on-matrices} Grant Assumptions
\ref{assassumption-on-noise}--\ref{assassumption-on-kernel}. Then, for
any $n\ge3$ and any $r\ge1$,
\[
\mathbb{P}^{(n)}_{X,\varepsilon} \Bigl\{\sup_{E\in\mathcal
E_*}
\bigl[ \bigl|\eta _{n,t}(E) \bigr|+ \bigl|\xi_{n,t}(E) \bigr| \bigr]\ge
C_1(n)\|F\|_\infty+C_2(n) \Bigr\}\le(8+
\Upsilon)n^{-4r}.
\]
The expressions for $ C_1(n) $ and $ C_2(n) $ are given in Section~\ref{secconstants}.
\end{lemma}
%
%
\begin{lemma}
\label{leminegality-for-sup-norm} Grant Assumptions
\ref{assassumption-on-noise}--\ref{assassumption-on-kernel}. Then, for
any $ n\ge n_0\vee n_1 $,
\[
\sup_{\theta^*\in\mathbb{S}^{1}}\sup_{f\in\mathbb{F}(\beta
_0,M)}\mathbb{P}^{(n)}_{F}
\bigl\{ \widehat{F}_\infty\notin \bigl[\| F\|_\infty,
3M+4C_5(n) \bigr] \bigr\}\le(8+\Upsilon)n^{-4r}.
\]
The numbers $ n_0, n_1 $ are defined
in~(\ref{eqrestriction-on-sample-size}) and $ C_5(n) $ is defined in
Section~\ref{secconstants}.
\end{lemma}

\subsubsection{\texorpdfstring{Proof of Theorem \protect\ref{thlocal-oracle-inequality}}
{Proof of Theorem 1}}
In view of Jensen's inequality, an upper bound for $ \mathcal
R_{r,t}^{(n)}
$,  $r\ge2 $, will suffice to complete the proof.

Let $h^*\in\mathcal H_n$ be a bandwidth such that $ 2h^*\le
h^*_{f}<4h^* $. Introduce the following random events:
\[
\mathcal A= \bigl\{R^{(1)}_{t}\bigl(\theta^*,h^*
\bigr)+R^{(2)}_{t}\bigl(h^*\bigr)=0 \bigr\},\qquad \mathcal
B= \bigl\{\widehat{F}_\infty\in \bigl[\|F\|_\infty,
3M+4C_5(n) \bigr] \bigr\}
\]
and let $ \widebar{\mathcal A} $ and $ \widebar{\mathcal B} $ denote
the events complimentary to $ \mathcal A $ and $ \mathcal B $,
respectively. The proof is split into three steps.

\subsection*{Risk computation under $ \mathcal A \cap\mathcal B $}
First, the following inclusion holds:
%
\begin{equation}
\label{eq00proof-of-theorem-local} \mathcal A\subseteq \bigl\{\hat{h}\ge h^* \bigr\}.
\end{equation}
Indeed, the definition of the couple $(\hat{\theta},\hat{h})$ yields
\begin{eqnarray*}
\mathrm{1}_{\mathcal A} \operatorname{TH}\bigl(h^*\bigr)&=&
\mathrm{1}_{\mathcal A} \bigl\{ R^{(1)}_{t}\bigl(\theta
^*,h^*\bigr)+R^{(2)}_{t}\bigl(h^*\bigr)+\operatorname{TH}
\bigl(h^*\bigr) \bigr\}
\\
&\ge& \mathrm{1}_{\mathcal A} \bigl\{R^{(1)}_{t}(
\hat{\theta},\hat {h})+R^{(2)}_{t}(\hat{h})+
\operatorname{TH}(\hat{h}) \bigr\} \ge \mathrm{1}_{\mathcal A}
\operatorname{TH}(\hat{h}).
\end{eqnarray*}
It remains to note that the mapping $\eta\mapsto\operatorname
{TH}(\eta)$ is decreasing so
inclusion~(\ref{eq00proof-of-theorem-local}) follows. Next, the
triangle inequality yields
%
\begin{eqnarray}
\label{eq0proof-of-theorem-local} \bigl\llvert \widehat F_{(\hat{\theta},\hat{h})}(t)-F(t)\bigr\rrvert &\le&
\bigl\llvert \widehat F_{(\theta^*,h^*)}(t)-F(t)\bigr\rrvert + \bigl\llvert
\widehat F_{(\hat
{\theta },\hat{h})}(t)- \widehat F_{(\hat{\theta},h^*)}(t)\bigr\rrvert
\nonumber
\\
&&{}+ \bigl\llvert \widehat F_{(\theta^*,h^*)(\hat{\theta},h^*)}(t)-\widehat
F_{(\hat{\theta},h^*)}(t) \bigr\rrvert
\\
&&{}+\bigl\llvert \widehat F_{(\theta^*,h^*)(\hat{\theta
},h^*)}(t)-\widehat
F_{(\theta^*,h^*)}(t) \bigr\rrvert.\nonumber
\end{eqnarray}
\begin{longlist}[$2^0$.]
\item[$1^0$.]  We have in view of (\ref{eq00proof-of-theorem-local})
    and
    the definition of $R^{(2)}_{t}$ that
%
\begin{equation}
\label{eq1proof-of-theorem-local} \mathrm{1}_{\mathcal A}\bigl\llvert \widehat
F_{(\hat{\theta},\hat
{h})}(t)- \widehat F_{(\hat{\theta},h^*)}(t)\bigr\rrvert \le
\mathrm{1}_{\mathcal A} \bigl[R^{(2)}_{t}(\hat{h})+
\operatorname{TH}\bigl(h^*\bigr) \bigr].
\end{equation}
The definition of $R^{(1)}_{t}(\cdot,\cdot)$ implies that
%
\begin{eqnarray}
\label{eq01proof-of-theorem-local} \mathrm{1}_{\mathcal A}\bigl\llvert \widehat
F_{(\theta^*,h^*)(\hat
{\theta},h^*)}(t)-\widehat F_{(\hat{\theta},h^*)}(t) \bigr\rrvert &\le& \mathrm{1}_{\mathcal A} \bigl[R^{(1)}_{t}\bigl(\theta^*,h^*
\bigr)+\operatorname {TH}\bigl(h^*\bigr) \bigr]
\nonumber
\nonumber\\[-8pt]\\[-8pt]
&=&\mathrm{1}_{\mathcal A}\operatorname{TH}\bigl(h^*\bigr).\nonumber
\end{eqnarray}
Note that $E_{(\theta, h)(\nu,h)}=\pm E_{(\nu, h)(\theta,h)}$, for any
$\theta,\nu$ and $h$. Hence,
\[
\widehat F_{(\theta^*,h^*)(\hat{\theta},h^*)}(\cdot)\equiv \widehat F_{(\hat
{\theta},h^*)(\theta^*,h^*)}(\cdot),
\]
because $\mathcal{K}$ is symmetric. The latter observation, inclusion
(\ref{eq00proof-of-theorem-local}) and the definition of $R^{(1)}_{t}$
yield
%
\begin{eqnarray}
\label{eq2proof-of-theorem-local} \mathrm{1}_{\mathcal A}\bigl\llvert \widehat
F_{(\theta^*,h^*)(\hat
{\theta},h^*)}(t)- \widehat F_{(\theta^*,h^*)}(t) \bigr\rrvert &=&
\mathrm{1}_{\mathcal A}\bigl\llvert \widehat F_{(\hat{\theta
},h^*)(\theta
^*,h^*)}(t)-\widehat
F_{(\theta^*,h^*)}(t) \bigr\rrvert
\nonumber\\[-8pt]\\[-8pt]
&\le&\mathrm{1}_{\mathcal A} \bigl[R^{(1)}_{t}(
\hat{\theta },\hat {h})+\operatorname{TH}\bigl(h^*\bigr) \bigr].\nonumber
\end{eqnarray}
From (\ref{eq0proof-of-theorem-local}),
(\ref{eq1proof-of-theorem-local}), (\ref{eq01proof-of-theorem-local})
and (\ref{eq2proof-of-theorem-local}), we obtain that
\begin{eqnarray*}
\mathrm{1}_{\mathcal A}\bigl\llvert \widehat F_{(\hat{\theta},\hat
{h})}(t)-F(t)\bigr
\rrvert &\le& \mathrm{1}_{\mathcal A} \bigl[R^{(1)}_{t}(
\hat{\theta },\hat{h})+ R^{(2)}_{t}(\hat{h})
\bigr]+3\operatorname{TH}\bigl(h^*\bigr)
\\
&&{} +\bigl\llvert \widehat F_{(\theta
^*,h^*)}(t)-F(t)\bigr\rrvert. 
\end{eqnarray*}
In addition, the definition of $ (\hat{\theta},\hat{h}) $ guarantees
that
\begin{eqnarray*}
R^{(1)}_{t}(\hat{\theta},\hat{h})+
R^{(2)}_{t}(\hat{h}) &\le& R^{(1)}_{t}(
\hat{\theta},\hat{h})+ R^{(2)}_{t}(\hat h)+
\operatorname{TH} (\hat{h})
\\
&\le& R^{(1)}_{t}\bigl(\theta^*,h^*\bigr)+
R^{(2)}_{t}\bigl(h^*\bigr)+\operatorname{TH}\bigl(h^*\bigr).
\end{eqnarray*}
We then obtain
%
\begin{eqnarray}
\label{eq02proof-of-theorem-local} \mathrm{1}_{\mathcal A}\bigl\llvert \widehat
F_{(\hat{\theta},\hat{h})}(t)-F(t)\bigr\rrvert &\le& 4\operatorname{TH}\bigl(h^*\bigr)+
\bigl\llvert \widehat F_{(\theta
^*,h^*)}(t)-F(t)\bigr\rrvert.
\end{eqnarray}
Note also that, for any $\eta\in\mathcal H_n$,
\begin{eqnarray*}
\mathrm{1}_{\mathcal B}\operatorname{TH}(\eta) &\le& 2 \bigl[\|\mathcal{K}
\|^{2}_{\infty}\sqrt{\ln(n)}+(3M+4C_5)
C_1(n)+C_2(n) \bigr](\eta n)^{-1/2}
\\
&\le& C_6\sqrt{(\eta n)^{-1}\ln(n)},
\end{eqnarray*}
where $C_6=2\|\mathcal{K}\|^{2}_{\infty}+2(3M+4C_5) C_1+2C_2$ and $C_1,
C_2$ and $C_5$ are defined in~(\ref{eqasympt-bound-for-quantities}).
Because $\operatorname{TH}(h^*)\le\operatorname{TH}(h_f^*/4)$, this
bound and (\ref{eq02proof-of-theorem-local}) yield
%
\begin{equation}
\label{eq020proof-of-theorem-local} \mathrm{1}_{\mathcal A\cap\mathcal B}\bigl\llvert \widehat
F_{(\hat
{\theta
},\hat{h})}(t)-F(t)\bigr\rrvert \le8 C_6\sqrt{
\frac{\ln(n)}{nh_f^*} }+\bigl\llvert \widehat F_{(\theta^*,h^*)}(t)-F(t)\bigr\rrvert.
\end{equation}

\item[$2^0$.]  Note that
    $E_{(\theta,h)(\nu,h)},E_{(\theta,h)}\in\mathcal
    E_*$, for any $\theta,\nu\in\mathbb{S}^{1}$, $h\in[h_{\min},1]$.
    Set
\[
\widehat{F}(E,t)=\frac{\det(E)}{n} \sum_{i=1}^n
K \bigl( E(X_i-t) \bigr)g^{-1}(X_i)
Y_i,\qquad E\in\mathcal E_*.
\]
The following ``approximation${}+{}$stochastic part'' decomposition of
$\widehat{F}(E,t)$ will be useful in the sequel:
%
\begin{eqnarray}
\label{eq03proof-of-theorem-local} \widehat{F}(E,t)&=&\det(E) \int K \bigl(E(x-t) \bigr)F(x)
\,\mathrm{d}x
\nonumber\\[-8pt]\\[-8pt]
&&{}+\sqrt{n^{-1}\det(E)} \bigl[\eta_{n,t}(E)+\xi_{n,t}(E) \bigr],\nonumber
\end{eqnarray}
where $\eta_{n,t}(E)$ and $\xi_{n,t}(E)$ are defined before the
statement of Lemma \ref{lemgauss-on-matrices}. Hence
\begin{eqnarray*}
&& \bigl\llvert \widehat F_{(\theta^*,h^*)}(t)-F(t)\bigr\rrvert
\\
&&\qquad \le
\bigl\llvert S_{(\theta^*, h^*)}(t) - F(t)\bigr\rrvert
+\sqrt{n^{-1}\det (E_{(\theta^*,h^*)} ) } \bigl|\eta _{n,t}
(E_{(\theta^*,h^*)} )+ \xi_{n,t} (E_{(\theta^*,h^*)} ) \bigr|.
\end{eqnarray*}
Taking into account that $\det (E_{(\theta^*,h^*)}
)=(h^*)^{-1}\le4(h^*_f)^{-1}$ in view of the definition of $h^*$ and
using the third assertion of Lemma \ref{lemboundsforbias}, we obtain
\begin{eqnarray*}
\label{eq04proof-of-theorem-local} &&\bigl\llvert \widehat F_{(\theta^*,h^*)}(t)-F(t)\bigr\rrvert
\\
&&\qquad \le\bigl(nh^*_f\bigr)^{-1/2} \bigl[\sqrt{\ln(n)}\|
\mathcal{K}\| _{\infty
}+2 \bigl|\eta_{n,t} (E_{(\theta^*,h^*)} )+
\xi_{n,t} (E_{(\theta^*,h^*)} ) \bigr| \bigr].
\end{eqnarray*}
Applying the Rosenthal inequality to $\eta_{n,t} (E_{(\theta ^*,h^*)}
)+\xi_{n,t} (E_{(\theta^*,h^*)} )$ which is a sum of centered
independent random variables, from (\ref{eq020proof-of-theorem-local})
we obtain
%
\begin{eqnarray}
\label{eq08proof-of-theorem-local} && \bigl\{\mathbb{E}_F^{(n)}\bigl\llvert
\widehat F_{(\hat{\theta
},\hat{h})}(t)-F(t)\bigr\rrvert ^r\mathrm{1}_{\mathcal A\cap\mathcal B}
\bigr\} ^{1/r} \le\tilde{c}_0\sqrt{
\bigl(nh_f^*\bigr)^{-1}\ln(n)},
\end{eqnarray}
where $\tilde{c}_0$ is independent of $F$ and $n$.
\end{longlist}

\subsection*{Risk computation under $ \widebar{\mathcal B} $}
Because $ f\in\mathbb{F}(\beta_0,M) $ and $ nh_{\min}>1 $, we have
\[
\bigl|\widehat F_{(\hat{\theta},\hat{h})}(t)-F(t) \bigr|\le n \biggl\{ M \bigl(1+
\underline{g}^{-1}\|K\|_\infty \bigr)+\underline{g}^{-1}
\| K \| _\infty n^{-1} \sum_{i=1}|
\varepsilon_i| \biggr\}.
\]
Hence, in view of the Rosenthal inequality, we obtain
%
\begin{equation}
\label{eq08-1proof-of-theorem-local} \bigl[\mathbb{E}_F^{(n)}\bigl\llvert
\widehat F_{(\hat{\theta
},\hat{h})}(t)-F(t)\bigr\rrvert ^{2r}
\bigr]^{1/(2r)} \le\tilde{c}_1 n,
\end{equation}
where $\tilde{c}_1$ is independent of $F$ and $n$.

The use of the Cauchy--Schwarz inequality together with the statement
of Lemma \ref{leminegality-for-sup-norm} leads to the following bound:
%
\begin{eqnarray}
\label{eq8proof-of-theorem-local} \bigl\{\mathbb{E}^{(n)}_{F} \bigl\llvert
\widehat F_{(\hat{\theta},\hat{h})}(t)-F(t)\bigr\rrvert ^r\mathrm{1}_{\widebar{\mathcal B}}
\bigr\} ^{1/r} &\le&\tilde{c}_1 n \bigl[
\mathbb{P}^{(n)}_F(\widebar{\mathcal B})
\bigr]^{1/(2r)}
\nonumber\\[-8pt]\\[-8pt]
&\le& \tilde{c}_1(8+\Upsilon)^{1/(2r)}n^{-1}.\nonumber
\end{eqnarray}
%

\subsection*{Risk computation under $ \widebar{\mathcal A}\cap
\mathcal B $} We note that
\[
\mathbb{P}^{(n)}_{F} \{\widebar{\mathcal A}\cap\mathcal B
\}\le \mathbb{P}^{(n)}_{F} \bigl\{R^{(1)}_{t}
\bigl(\theta^*,h^*\bigr)>0, \mathcal B \bigr\} +\mathbb{P}^{(n)}_{F}
\bigl\{R^{(2)}_{t}\bigl(h^*\bigr)>0, \mathcal B \bigr\}.
\]

\begin{longlist}[$2^0$.]
\item[$1^0$.]  First, let us bound from above $\mathbb{P}^{(n)}_{F}
    \{R^{(1)}_{t}(\theta^*,h^*)>0, \mathcal B \}$. We have
\[
\bigl\{R^{(1)}_{t}\bigl(\theta^*,h^*\bigr)>0 \bigr\}=
\bigcup_{\eta\in\mathcal H_n\dvtx  \eta\le h^*} \Bigl\{\sup_{\nu\in\mathbb{S}^{1}} \bigl|
\widehat F_{(\theta^*,
\eta
)(\nu, \eta
)}(t) - \widehat F_{(\nu,\eta)}(t) \bigr|> \operatorname{TH}
(\eta ) \Bigr\}
\]
and, therefore,
%
\begin{eqnarray}\label{eq9proof-of-theorem-local}
\qquad && \mathbb{P}^{(n)}_{F} \bigl\{
R^{(1)}_{t}\bigl(\theta^*,h^*\bigr)>0, \mathcal B \bigr\}\nonumber
\\
&&\qquad \le\sum_{k\dvtx  2^{-1}h_{\min}\le2^{-k}\le h^*}
\mathbb{P}_F^{(n)} \Bigl\{\sup
_{\nu\in\mathbb{S}^{1}} \bigl|\widehat F_{(\theta^*, 2^{-k})(\nu, 2^{-k})}(t) - \widehat
F_{(\nu,2^{-k})}(t) \bigr|
\\
&&\hspace*{226pt} > \operatorname{TH} \bigl(2^{-k} \bigr), \mathcal B \Bigr\}.\nonumber
\end{eqnarray}
Thus, denoting by $\varsigma_n=\sup_{E\in\mathcal E_*} [ |\eta_{n,t}(E)
|+ |\xi_{n,t}(E) | ]$ and using (\ref{eq03proof-of-theorem-local})
together with the first assertion of Lemma \ref{lemboundsforbias}, we
obtain, for any $k\dvtx 2^{-k}\le h^*$,
%
\begin{eqnarray}
\label{eq10proof-of-theorem-local} &&\sup_{\nu\in\mathbb{S}^{1}}\bigl\llvert \widehat
F_{(\theta^*,2^{-k})(\nu,2^{-k})}(t)-\widehat F_{(\nu,2^{-k})}(t)\bigr\rrvert
\nonumber
\\
&&\qquad \le2\bigl(h^*_f\bigr)^{-1/2}\| \mathcal{K}
\|^{2}_{\infty}\sqrt{n^{-1}\ln (n)}+2\sqrt
{2^k}n^{-1/2} \varsigma_n
\\
&&\qquad \le2 \| \mathcal{K}\|^{2}_{\infty}\sqrt{2^kn^{-1}
\ln(n)}+2\sqrt {2^{k}n^{-1}} \varsigma_n.\nonumber
\end{eqnarray}
Here we have also used that $2^{-1}h^*_f\ge2^{-k}$. Note also that
%
\begin{eqnarray}
\label{eq010proof-of-theorem-local} &&\mathrm{1}_{\mathcal B}\operatorname{TH}(\eta)\ge2\|
\mathcal{K}\| ^{2}_{\infty} \sqrt{\frac{\ln(n)}{\eta n}}+
\frac{2}{\sqrt{\eta n}} \bigl( C_1(n) \|F\|_\infty+C_2(n)
\bigr)
\end{eqnarray}
and, therefore, we obtain from (\ref{eq10proof-of-theorem-local}), for
any $k$ satisfying $2^{-k}\le h^*$,
\begin{eqnarray*}
&& \mathbb{P}^{(n)}_{F} \Bigl\{\sup_{\nu\in\mathbb{S}^{1}}
\bigl|\widehat F_{(\theta^*, 2^{-k})(\nu, 2^{-k})}(t) - \widehat F_{(\nu,
2^{-k})}(t) \bigr|>
\operatorname{TH} \bigl(2^{-k} \bigr), \mathcal B \Bigr\}
\nonumber
\\
&&\qquad \le\mathbb{P}^{(n)}_{X,\varepsilon} \bigl\{
\varsigma_n\ge\|F\| _\infty C_1(n)+C_2(n)
\bigr\}\le(8+\Upsilon)n^{-4r},
\end{eqnarray*}
in view of Lemma \ref{lemgauss-on-matrices}. This bound and
(\ref{eq9proof-of-theorem-local}) yield
%
\begin{equation}
\label{eq11proof-of-theorem-local} \mathbb{P}^{(n)}_{F} \bigl
\{R^{(1)}_{t}\bigl(\theta^*,h^*\bigr)>0, \mathcal B \bigr\}
\le (8+\Upsilon) \log_2 (n) n^{-4r}.
\end{equation}

\item[$2^0$.]  Now, let us bound from above $\mathbb{P}^{(n)}_{F}
    \{R^{(2)}_{t}(h^*)>0,  \mathcal B \}$. We have
\[
\bigl\{R^{(2)}_{t}\bigl(h^*\bigr)>0 \bigr\}=\bigcup
_{\eta\in\mathcal H_n\dvtx  \eta
\le h^*} \Bigl\{\sup_{\theta\in\mathbb{S}^{1}}\bigl
\llvert \widehat F_{(\theta,h^*)}(t)- \widehat F_{(\theta, \eta)}(t)\bigr\rrvert >
\operatorname{TH} \bigl(2^{-k} \bigr) \Bigr\}
\]
and, hence,
%
\begin{eqnarray}
\label{eq12proof-of-theorem-local}
&& \mathbb{P}^{(n)}_{F} \bigl\{R^{(2)}_{t}\bigl(h^*\bigr)>0, \mathcal B \bigr\}\nonumber
\\
&&\qquad \le\sum_{k\dvtx  2^{-1}h_{\min}\le2^{-k}\le h^*}\mathbb{P}^{(n)}_{F}
\Bigl\{\sup_{\theta\in\mathbb{S}^{1}}\bigl\llvert \widehat F_{(\theta,h^*)}(t)-
\widehat F_{(\theta,2^{-k})}(t)\bigr\rrvert
\\
&&\hspace*{192pt} > \operatorname{TH}\bigl(2^{-k} \bigr), \mathcal B \Bigr\}.\nonumber
\end{eqnarray}
Similar to estimate (\ref{eq10proof-of-theorem-local}), with the use
of (\ref{eq03proof-of-theorem-local}) and the second assertion of
Lemma~\ref{lemboundsforbias}, we obtain, for any $k$ satisfying
$2^{-k}\le h^*$, that
\[
\sup_{\theta\in\mathbb{S}^{1}}\bigl\llvert \widehat F_{(\theta,h^*)}(t)-
\widehat F_{(\theta,
2^{-k})}(t)\bigr\rrvert \le2 \| \mathcal{K}\|^{2}_{\infty}
\sqrt{2^kn^{-1}\ln(n)}+2\sqrt {2^{k}n^{-1}}
\varsigma_n.
\]
For any $k$ satisfying $2^{-k}\le h^*$, bound
(\ref{eq010proof-of-theorem-local}) and Lemma~\ref{lemgauss-on-matrices} yield
\begin{eqnarray*}
&&\mathbb{P}^{(n)}_{F} \Bigl\{\sup_{\theta\in\mathbb{S}^{1}}
\bigl\llvert \widehat F_{(\theta,h^*)}(t)- \widehat F_{(\theta, 2^{-k})}(t)\bigr
\rrvert > \operatorname {TH} \bigl(2^{-k} \bigr), \mathcal B \Bigr\}
\nonumber
\\
&&\qquad \le\mathbb{P}^{(n)}_{X,\varepsilon} \bigl\{
\varsigma_n\ge\|F\| _\infty C_1(n)+C_2(n)
\bigr\}\le(8+\Upsilon) n^{-4r}.
\end{eqnarray*}
Together with (\ref{eq12proof-of-theorem-local}), the latter bound
gives
%
\begin{equation}
\label{eq15proof-of-theorem-local} \mathbb{P}^{(n)}_{F} \bigl
\{R^{(2)}_{t}\bigl(h^*\bigr)>0, \mathcal B \bigr\}\le (8+
\Upsilon ) \log_2 (n) n^{-4r}.
\end{equation}
Thus, we obtain from (\ref{eq11proof-of-theorem-local}) and
(\ref{eq15proof-of-theorem-local}) that
\[
\mathbb{P}^{(n)}_{F} (\widebar{\mathcal A}\cap\mathcal B )
\le 2(8+\Upsilon) \log_2 (n) n^{-4r}.
\]
Subsequently, this bound and (\ref{eq08-1proof-of-theorem-local}) yield
%
\begin{equation}
\label{eq16proof-of-theorem-local}
\qquad\quad \bigl\{\mathbb{E}^{(n)}_{F}\bigl\llvert
\widehat F_{(\hat{\theta
},\hat{h})}(t)-F(t)\bigr\rrvert ^r\mathrm{1}_{\widebar{\mathcal A}\cap\mathcal
B}\bigr\}^{1/r} \le\tilde{c}_1 n \bigl[
\mathbb{P}^{(n)}_{F}(\widebar{\mathcal A}\cap \mathcal B)
\bigr]^{1/(2r)}\le\tilde{c}_2n^{-1/2},
\end{equation}
where $\tilde{c}_2$ is independent of $F$ and $n$.

The assertion of the theorem follows from
(\ref{eq08proof-of-theorem-local}), (\ref{eq8proof-of-theorem-local})
and (\ref{eq16proof-of-theorem-local}).
\end{longlist}

\subsection{\texorpdfstring{Proof of Theorem \protect\ref{thpointwise-adaptation}}
{Proof of Theorem 3}}
Using the standard computation of the bias of kernel estimators, under
Assumptions \ref{assassumption-on-kernel} and
\ref{ass2assumption-on-kernel}, we get, for any $f\in\mathbb
{H}(\beta,L)$ and any $z\in\mathbb{R}$,
\[
\Delta_{\mathcal{K},f}(h,z) 
\le \frac{ L h^\beta2^{-\beta} \|\mathcal{K}\|_\infty}{ (1+\beta
)m_{\beta}! } \le
\|\mathcal{K}\|_\infty Lh^\beta.
\]
Since the right-hand side of the latter inequality is independent of
$z$, we have $ \Delta^*_{\mathcal{K},f}(h,z) \le\|\mathcal{K}\|_\infty
Lh^\beta$.  This\vspace*{1pt} implies $h^*_{\mathcal{K},f}(z)\ge
(L^{-2}n^{-1}\ln(n) )^{1/(2\beta+1)}$, for any~$z\in\mathbb{R}$, so the
assertion of the theorem follows from Theorem
\ref{thlocal-oracle-inequality}.

\subsection{\texorpdfstring{Proof of Theorem \protect\ref{thpointwise-adaptation-lower}}
{Proof of Theorem 4}}
We start this section with an auxiliary result used in the proof of the
second assertion of the theorem. It was established in
\citet{KLP2008}, Corollary 2 of Proposition 5, and, for
convenience, we formulate it as Lemma \ref{lemKLP-result} below.

\subsubsection{Auxiliary result}
The result cited below concerns a lower bound for estimators of an
arbitrary mapping in the framework of an abstract statistical model. We
do not present it in full generality and below a version reduced to the
estimation at a given point is provided.

Let $ \mathcal F $ be a nonempty class of functions; and let $ F\dvtx
\mathbb{R}^d\to\mathbb{R}$ be an unknown function from model defined in
(\ref{model-regression})--(\ref{single-index}). The aim is to estimate
the functional $ F(t)$,  $t \in[-1/2, 1/2]^d$.

Introduce the following notation. For any given $F, G\in\mathcal{F}$,
set
\[
Z(F,G)=\prod_{i=1}^n
\biggl[\frac{p (Y_i-F(X_i) )}{p
(Y_i-G(X_i) )} \biggr].
\]

%
\begin{lemma}
\label{lemKLP-result} Assume that, for any sufficiently large $ n \ge
1$, there exist a positive integer $N_{n}$,  $c>1$ and functions $
F_0,\ldots, F_{N_{n}} \in\mathcal F $ such that
%
\begin{eqnarray}
\label{eqass1-klp-lemma} \bigl|F_j(t)-F_0(t)\bigr| &=&\lambda_n\qquad \forall j=1, \ldots, N_{n},
\\
\label{eqass2-klp-lemma} \mathbb{E}^{(n)}_{F_0} \Biggl(
\frac{1}{N_{n}}\sum_{j=1}^{N_{n}}Z
(F_j,F_0 ) \Biggr)^{2}&\le& c. 
\end{eqnarray}
Then, for $ r\ge1 $ and any $ t \in[-1/2, 1/2]^d$,
\[
\label{lowerboundGWN} \inf_{\widetilde F} \sup_{F \in\mathcal F}
\bigl(\mathbb{E}^{(n)}_{F} \bigl| \widetilde{F} (t) - F (t)
\bigr|^r \bigr)^{1/r} \ge 
\frac{1}{2} \biggl[1- \sqrt{\frac{c-1}{c+3}} \biggr]
\lambda_n.
\]
\end{lemma}

\subsubsection{\texorpdfstring{Proof of Theorem \protect\ref{thpointwise-adaptation-lower}}
{Proof of Theorem 4}}
The proof is based on the construction of $F_0,\ldots,F_{N_{n}}$
satisfying conditions
(\ref{eqass1-klp-lemma})--(\ref{eqass2-klp-lemma}) of Lemma
\ref{lemKLP-result}.
\begin{longlist}[$2^0$.]
\item[$1^0$.]   First, we construct $ F_0,\ldots,F_{N_{n}} $ and verify
    (\ref{eqass1-klp-lemma}). Let $ w\dvtx \mathbb{R}\to\mathbb{R}$ be
    a~function such that $ \operatorname{supp}(w) \subset(-1/2, 1/2)
$,  $w \in\mathbb{H}(\beta,1) $ and $w(0)\neq0$. Set $ h=
(\mathfrak{a}(L^2n)^{-1}\ln(n) )^{1/(2\beta+1)} $, where
$\mathfrak{a}>0$ will be chosen later, and define
%
\begin{equation}
\label{linkfunctionforhypothesis} f(z) = Lh^\beta w \bigl(zh^{-1} \bigr),\qquad z\in
\mathbb{R}.
\end{equation}
For $ b>0 $, put $ N_{n}= n^{b} $ assuming without loss of generality
that $ N_{n}$ is an integer. The value of $ b $ will be determined
later in order to satisfy (\ref{eqass2-klp-lemma}).

Let $ \{\vartheta_j,  j=1,\ldots, N_{n} \}\subset \mathbb{S}^{d-1} $ be
defined as follows:
\[
\vartheta_j= \bigl(\theta^{(1)}_j,
\theta^{(2)}_j,0,\ldots,0 \bigr)^{\top}, \qquad
\theta^{(1)}_j=\cos(j/N_{n}),
\theta^{(2)}_j=\sin(j/N_{n}).
\]
Finally, we set
%
\begin{equation}
\label{hypotheses} F_0 \equiv0\quad\mbox{and}\quad F_j(x) =
f \bigl(\vartheta_j^{\top} (x-t) \bigr),\qquad j = 1, \ldots,
N_{n}.
\end{equation}
Obviously, $ f $ defined by (\ref{linkfunctionforhypothesis}) belongs
to $ \mathbb{H}(\beta,L) $, so all $ F_i $ are in the class $ \mathcal
F= \mathbb{F}_d(\beta,L) $. Moreover, for any $ i=1,\ldots, N_{n}$,
\begin{eqnarray*}
\bigl|F_j(t) - F_0(t)\bigr| &=&
\bigl|w(0)\bigr| L^{1/(2\beta+1)} \bigl(\mathfrak{a}n^{-1}\ln(n)
\bigr)^{\beta/(2\beta+1)}
\\
&=& \bigl|w(0)\bigr|\mathfrak{a}^{\beta/(2\beta+1)}\psi_n(\beta,L).
\end{eqnarray*}
We see that (\ref{eqass1-klp-lemma}) holds with $ \lambda_n = |w(0)|
\mathfrak{a}^{\beta/(2\beta+1)}\psi_n(\beta,L) $.

\item[$2^0$.]   It is noteworthy that
\begin{eqnarray*}
&&\mathbb{E}^{(n)}_{F_0} \Biggl[\frac{1}{N_{n}}\sum
_{j=1}^{N_{n}}Z (F_j,F_0
) \Biggr]^{2}
\\
&&\qquad= \frac{1}{N^2_n}\sum_{j=1}^{N_{n}}
\mathbb {E}^{(n)}_{F_0} \bigl[Z^{2}
(F_j,F_0 ) \bigr] + \frac{1}{N^2_n}\mathop{\sum
_{j,k=1,}}_{j\neq k}^{N_{n}}
\mathbb{E}^{(n)}_{F_0} \bigl[Z (F_j,F_0
)Z (F_k,F_0 ) \bigr].
\end{eqnarray*}
It follows that
\begin{eqnarray*}
\mathbb{E}^{(n)}_{F_0} \bigl[Z^{2}
(F_j,F_0 ) \bigr] &=& \biggl\{\int_{\mathbb{R}^{d+1}}
\frac{p^{2} (y-F_j(x)
)}{p(y)}g(x)\,\mathrm{d}x\,\mathrm{d}y \biggr\}^n,
\\
\mathbb{E}^{(n)}_{F_0} \bigl[Z (F_j,F_0
)Z (F_k,F_0 ) \bigr] &=& \biggl\{\int_{\mathbb{R}^{d+1}}
\frac{p (y-F_j(x) )p
(y-F_k(x) )}{p(y)}g(x)\,\mathrm{d}x\,\mathrm{d}y \biggr\}^n.
\end{eqnarray*}
Because $ \lim_{n\to\infty} \sup_{j=1,\ldots,N_{n}}\|F_j\|_\infty=0$,
we have in view of Assumptions~\ref{ass2on-noise}
and~\ref{ass2on-design}, for all $n$ large enough,
%
\begin{eqnarray}
&& \int_{\mathbb{R}^{d+1}} \biggl[\frac{p^{2} (y-F_j(x)
)}{p(y)} \biggr]g(x)\,\mathrm{d}x
\,\mathrm{d}y\nonumber
\\
&&\qquad \le 1+\mathfrak{Q}\int_{\mathbb{R}^d}F^{2}_j(x)g(x)
\,\mathrm{d}x \label{eq01regression-lower}
\\
&&\qquad \le 1+\mathfrak{Q}\mathfrak{g}\int
_{\mathbb{R}^d}F^{2}_j(x)
\bigl(1+|x|_2^{\varpi} \bigr)^{-1}\,\mathrm{d}x,\nonumber
\\
&& \int_{\mathbb{R}^{d+1}} \biggl[\frac{p (y-F_j(x) )p (y-F_k(x) )}{p(y)} \biggr]g(x)
\,\mathrm{d}x\,\mathrm{d}y\nonumber
\\
&&\qquad \le 1+ \mathfrak{Q}\int_{\mathbb{R}^d}
\bigl|F_j(x)F_k(x) \bigr|g(x)\,\mathrm{d}x\label{eq02regression-lower}
\\
&&\qquad \le 1+ \mathfrak{Q}\mathfrak{g}\int
_{\mathbb{R}^d} \bigl|F_j(x)F_k(x) \bigr|
\bigl(1+|x|_2^{\varpi} \bigr)^{-1}\,\mathrm{d}x.\nonumber
\end{eqnarray}
Set ${\theta_{j}}_{\bot}= (-\sin(j/N_{n}),\cos(j/N_{n}) )^{\top}$ and
${\vartheta_j}_{\bot}= ({\theta_{j}}_{\bot}^{\top},0,\ldots,0
)^{\top}\in\mathbb{S}^{d-1}$. Denote for all $ j=1,\ldots,N_{n}$ by
$\Theta_j^{\top}$ the\vspace*{-1pt} orthogonal matrix $
(\vartheta_j,{\vartheta_j}_{\bot},\mathbf{e}_3,\ldots,\mathbf {e}_d )$,
where $\mathbf{e}_s$, $s=3,\ldots, d$, are the canonical basis vectors
in $\mathbb{R}^d$. Integration by substitution with $\Theta_j x=v$
gives
\begin{eqnarray*}
\int_{\mathbb{R}^d}F^{2}_j(x)
\bigl(1+|x|_2^{\varpi} \bigr)^{-1}\,\mathrm{d}x &=&
L^2h^{2\beta} \int_{\mathbb{R}^d}w^{2}
\bigl[ h^{-1} \bigl(v_1 - \vartheta _j^{\top
}
t\bigr) \bigr] \bigl(1+| v|_2^{\varpi} \bigr)^{-1}
\,\mathrm{d}v
\nonumber
\\
&\le& C_{\varpi}L^2\|w\|^2_2
h^{2\beta+1}=\mathfrak{a}C_{\varpi} \|w\| ^2_2n^{-1}
\ln(n),
\end{eqnarray*}
where we have denoted $C_{\varpi}=\int_{\mathbb{R}^{d-1}}
(1+|\mathrm{v}|_2^{\varpi} )^{-1}\,\mathrm{d} \mathrm{v}$ and
$\mathrm{v}=(v_2,\ldots,v_d)^{\top}$. For $ n $ sufficiently large,
this bound, together with (\ref{eq01regression-lower}), leads to
%
\begin{equation}
\label{eq2regression-lower} \sup_{j=1,\ldots,N_{n}}\mathbb{E}^{(n)}_{F_0}
\bigl[Z^{2} (F_j,F_0 ) \bigr] \le
n^{\mathfrak{a} \mathfrak{Q}\mathfrak{g}C_{\varpi} \|w\|^2_2}.
\end{equation}

For any $j\neq k$, set $\Theta_{j,k}^{\top} =
(\vartheta_j,\vartheta_k,\mathbf{e}_3,\ldots,\mathbf{e}_d )$. By
changing of variables with $\Theta_{j,k}x=v$, we have
\begin{eqnarray*}
&&\int_{\mathbb{R}^d} \bigl|F_j(x)F_k(x) \bigr|
\bigl(1+|x|_2^{\varpi
} \bigr)^{-1}\,\mathrm{d}x
\\
&&\qquad = \bigl|\det(\Theta_{j,k})\bigr|^{-1} L^2
h^{2\beta} \int_{\mathbb{R}^d} \frac{ | w [h^{-1}(v_1 -\vartheta
_j^{\top}t)
]w [h^{-1}(v_2 -\vartheta_k^{\top}t) ] |} {
1+ |\Theta^{-1}_{j,k}v |_2^{\varpi}}\,\mathrm{d}v
\\
&&\qquad \le \bigl|\det(\Theta_{j,k})\bigr|^{-1} c_{\varpi}L^2
h^{2\beta+2}\|w\|^2_1, 
\end{eqnarray*}
where $c_{\varpi}=\int_{\mathbb{R}^{d-2}} (1+|\mathrm{v}|_2^{\varpi }
)^{-1}\,\mathrm{d} \mathrm{v}$ and $ \mathrm{v}=(v_3,\ldots,v_d)^{\top}
$. Note that
\begin{eqnarray*}
\bigl|\det(\Theta_{j,k})\bigr| &=& \bigl|\cos(j/N_{n})
\sin(k/N_{n})-\cos (k/N_{n})\sin(j/N_{n}) \bigr|
\\
&=&\bigl|\sin \bigl((k-j)/N_{n} \bigr) \bigr|\ge\sin (1/N_{n} )>
(2N_{n})^{-1}
\end{eqnarray*}
for sufficiently large $n$.
We obtain
\[
\int_{\mathbb{R}^d} \bigl|F_j(x)F_k(x) \bigr|
\bigl(1+|x|_2^{\varpi
} \bigr)^{-1}\,\mathrm{d}x\le2
\mathfrak{a}c_{\varpi}\|w\|^2_1 n^{-1}
\ln (n)N_{n}h.
\]
Hence, choosing $b<1/(2\beta+1)$, we obtain, for all $n$ large enough,
that
%
\begin{equation}
\label{eq1proof-the-local-adaptive}
\qquad \sup_{j \ne k;  j,k = 1, \ldots, N_{n}} \int_{\mathbb
{R}^d}
\bigl|F_j(x)F_k(x) \bigr| \bigl(1+|x|_2^{\varpi}
\bigr)^{-1}\,\mathrm{d}x \le2\mathfrak{a}c_{\varpi}\|w
\|^2_1 n^{-1}.
\end{equation}
We have in view of (\ref{eq02regression-lower}) and
(\ref{eq1proof-the-local-adaptive})
%
\begin{eqnarray}
\label{eq001proof-the-local-adaptive} \sup_{j \ne k;  j,k = 1, \ldots, N_{n}} \mathbb{E}^{(n)}_{F_0}
\bigl\{ Z (F_j,F_0 )Z (F_k,F_0
) \bigr\} & \le& e^{2\mathfrak{a} \mathfrak{Q}\mathfrak{g}c_{\varpi}\|w\|^2_1}
\end{eqnarray}
and, hence, (\ref{eq2regression-lower}) and
(\ref{eq001proof-the-local-adaptive}) give
\[
\mathbb{E}^{(n)}_{F_0} \Biggl(\frac{1}{N_{n}}\sum
_{j=1}^{N_{n}}Z (F_j,F_0 )
\Biggr)^{2} \le n^{-b+\mathfrak{a} \mathfrak{Q}\mathfrak{g}C_{\varpi}\|w\|^2_2} +e^{2\mathfrak{a} \mathfrak{Q}\mathfrak{g}c_{\varpi}\|w\|^2_1}.
\]
Choosing $\mathfrak{a}=b (\mathfrak{Q}\mathfrak{g}C_{\varpi}\| w\|^2_2
)^{-1}$, we see that (\ref{eqass2-klp-lemma}) holds with the constant $c=1+
e^{2\mathfrak{a} \mathfrak{Q}\mathfrak{g}c_{\varpi}\|w\|^2_1}$. Since
$c$ appearing in (\ref{eqass2-klp-lemma}) is chosen
independently of $L$, the assertion of the theorem follows from
Lemma~\ref{lemKLP-result}.
\end{longlist}

\subsubsection{\texorpdfstring{Proof of Theorem \protect\ref{thglobal-adaptation}}
{Proof of Theorem 5}}

In the proof we exploit the ideas from \citet{LepMamSpok97}.
Moreover, our considerations are, to a great degree, based on the
technical result of Lemma \ref{lemlp-norm-of-bais} below. Its proof is
moved to the supplementary material [\citet{suppA}].

%
\begin{lemma}
\label{lemlp-norm-of-bais} Grant Assumptions
\ref{assassumption-on-kernel} and \ref{ass2assumption-on-kernel}. Then,
for any $\mathfrak{p}> 1$, $0<s\le\mathbf{b}$ and $\mathcal {Q}>0$, we
have
\[
\sup_{g\in\mathbb{N}_\mathfrak{p}(s,\mathcal{Q})} \bigl\llVert \Delta^*_{\mathcal{K},g}(h,\cdot)
\bigr\rrVert _\mathfrak{p}\le 2\tau _\mathfrak{p}\mathcal{Q}h^{s}
\|\mathcal{K}\|_\infty \bigl[2^{s \mathfrak{p}}-1 \bigr]^{-1/\mathfrak{p}}\qquad \forall h>0.
\]
Here $\tau_\mathfrak{p}$ is a dependent only of $ \mathfrak{p} $
constant from the $(\mathfrak{p},\mathfrak{p})$-strong maximal
inequality; see \citet{WheedenZygmund1977} for more details.
\end{lemma}

\begin{pf*}{Proof of Theorem \protect\ref{thglobal-adaptation}}
It is sufficient to prove the theorem in the case $r\ge p$ only.
Indeed, let us recall that the risk $\mathcal R_{r}^{(n)}(\cdot,
\cdot)$ is described by the $ L_r $ norm on $[-1/2,1/2]$, therefore,
\[
\mathcal R_{r}^{(n)}(\cdot, \cdot)\le\mathcal
R_{p}^{(n)}(\cdot, \cdot), r\le p.
\]
Hence, the case $r\le p$ can be reduced to the case $r=p$.

In view of Theorem \ref{thglobal-oracle-inequality}, in order to obtain
the assertion of the theorem, it suffices to bound from above $ \|
(nh^*_{\mathcal{K},f})^{-1}\ln(n) \|^{1/2}_{r/2}$.

Set\vspace*{-1pt} $\Gamma_0= \{y\in[-1/2,1/2]\dvtx h^*_{\mathcal
{K},f}(y)=1 \} $ and $\Gamma_k= \{y \in[-1/2,1/2]\dvtx\break
h^*_{\mathcal{K},f}(y) \in(2^{-k}, 2^{-k+1}] \cap [h_{\min},1 ] \}$,
for $ k=1,2,\ldots.  $ In what follows, the integration over the empty
set is supposed to be zero. We have
\[
\biggl\|\frac{\ln(n)}{nh^*_{\mathcal{K},f}} \biggr\|^{r/2}_{r/2} = \sum
_{k\ge1}\int_{\Gamma_k} \biggl(\frac{\ln(n)}{nh^*_{\mathcal{K},f}(y)}
\biggr)^{r/2}\,\mathrm{d}y+\int_{\Gamma_0} \biggl(
\frac
{\ln
(n)}{nh^*_{\mathcal{K},f}(y)} \biggr)^{r/2}\,\mathrm{d}y.
\]
For simplicity of notation, we denote by $ \bar{c}_i
$,  $i \ge1 $, constants independent of $n$, $f$ and $L$.

The definition of $\Gamma_0$ implies
%
\begin{equation}
\label{eq2proof-theorem-global-adaptation} \int_{\Gamma_0} \biggl(\frac{\ln(n)}{nh^*_{\mathcal{K},f}(y)}
\biggr)^{r/2}\,\mathrm{d}y\le\bar{c}_1
\bigl[n^{-1}\ln(n) \bigr]^{r/2}.
\end{equation}
We have in view of (\ref{eq1deforaclebandwidth}), for any $k\ge1$,
%
\begin{equation}
\label{eq3proof-theorem-global-adaptation} \Delta^*_{\mathcal{K},f} \bigl( h^*_{\mathcal{K},f}(y),y \bigr)=
\biggl[\frac{\| \mathcal{K}\|
^2_{\infty} \ln(n)}{n h^*_{\mathcal{K},f}(y)} \biggr]^{1/2}\qquad\forall y\in
\Gamma_k.
\end{equation}
Let $0\le q_k\le r$ be a sequence whose choice will be done later. We
obtain from~(\ref{eq3proof-theorem-global-adaptation})
%
\begin{eqnarray}
\label{eq4proof-theorem-global-adaptation}
&&\sum_{k\ge1}\int
_{\Gamma_k} \biggl(\frac{\ln(n)}{nh^*_{\mathcal{K},f}(y)} \biggr)^{r/2}
\,\mathrm{d}y\nonumber
\\
&&\qquad \le\bar{c}_2\sum_{k\ge1}
\biggl(\frac{ \ln(n)}{n2^{-k}} \biggr)^{(r-q_k)/2} \int_{\Gamma_k}
\bigl(\Delta^*_{\mathcal{K},f} \bigl(2^{1-k},y \bigr)
\bigr)^{q_k}\,\mathrm{d}y
\\
&&\qquad \le\bar{c}_2\sum
_{k\ge1} \biggl(\frac{ \ln(n)}{n2^{-k}} \biggr)^{(r-q_k)/2} \int
\bigl( \Delta^*_{\mathcal{K},f} \bigl(2^{1-k},y \bigr)
\bigr)^{q_k}\,\mathrm{d} y=:\Xi.\nonumber
\end{eqnarray}
To get the first inequality, we have used that $\Delta^*_{\mathcal
{K},f} (\cdot,y )$ is a monotonically increasing function.

The computation of the quantity on the right-hand side of
(\ref{eq4proof-theorem-global-adaptation}), including the choice of
$(q_k, k\ge1)$, will be done differently in dependence on $\beta, p$
and $r$.
\begin{longlist}[$2^0$.]
\item[$1^0$.]  \textit{Case} $(2\beta+1)p > r$. Put $ h^*=
    [L^{-2}n^{-1} \ln(n) ]^{1/(2\beta+1)} $ and choose $q_k=p$ if
    $2^{-k}\le h^*$ and $q_k=0$ if $2^{-k}> h^*$. By applying Lemma
    \ref{lemlp-norm-of-bais} with $\mathfrak{p}=p$, $s=\beta$ and
    $\mathcal{Q}=L$, we get
%
\begin{eqnarray}
\label{eq55proof-theorem-global-adaptation} 
\qquad\qquad \Xi &\le& \bar{c}_3L^{p}
\sum_{k\dvtx  2^{-k}\le h^*} \biggl(\frac{
\ln(n)}{n2^{-k}}
\biggr)^{(r-p)/2}2^{-k\beta p}+\bar{c}_4 \biggl(
\frac{\ln(n)}{nh^*} \biggr)^{r/2}
\nonumber\\[-8pt]\\[-8pt]
&\le& \bar{c}_5 \biggl[L^{p}
\bigl(n^{-1} \ln(n) \bigr)^{(r-p)/2}\sum
_{k\dvtx  2^{-k}\le h^*}2^{-k[\beta p-(r-p)/2 ]} + \biggl(\frac{\ln(n)}{nh^*}\biggr)^{r/2} \biggr].\nonumber
\end{eqnarray}
Because in the considered case $\beta p-\frac{r-p}{2}>0$, we obtain
\[
\Xi\le\bar{c}_6 \biggl[L^{ p}
\bigl(n^{-1} \ln(n) \bigr)^{(r-p)/2}\bigl(h^*\bigr)^{\beta p-(r- p)/2}
+ \biggl(\frac{
\ln(n)}{nh^*} \biggr)^{r/2} \biggr].
\]
It remains to note that $h^*$ is chosen by balancing two terms on the
right-hand side of the latter inequality. It yields
%
\begin{equation}
\label{eq5proof-theorem-global-adaptation} \Xi\le\bar{c}_7L^{r/(2\beta+1)}
\bigl(n^{-1}\ln (n) \bigr)^{(r\beta)/(2\beta+1)}.
\end{equation}
The argument in the case $(2\beta+1)p> r$ is completed with the use of
Theorem~\ref{thglobal-oracle-inequality},
(\ref{eq2proof-theorem-global-adaptation}) and
(\ref{eq5proof-theorem-global-adaptation}).

\item[$2^0$.]  \textit{Case} $(2\beta+1)p=r$. Put $h^*=1$ and choose
    $q_k=p$ for all $k\ge1$. Repeating the computations that led to
    (\ref{eq55proof-theorem-global-adaptation}), we get
%
\begin{equation}
\label{eq6proof-theorem-global-adaptation} \Xi\le\bar{c}_8\ln(n)L^{ p}
\bigl(n^{-1} \ln(n) \bigr)^{(r-p)/2}.
\end{equation}
Here we have used that $\beta p-(r-p)/2=0$ and that the summation
in (\ref{eq4proof-theorem-global-adaptation}) runs over $k$ such that
$2^{-k}\ge h_{\min}$, since otherwise $\Gamma_k=\varnothing$. It
remains to note that the equality $(2\beta+1)p=r$ is equivalent to
$p/r=1/(2\beta+1)$ and $(r-p)/2r=\beta/(2\beta+1)$. The assertion of
the theorem in the case $(2\beta+1)p=r$ follows now from
Theorem~\ref{thglobal-oracle-inequality},
(\ref{eq2proof-theorem-global-adaptation}) and
(\ref{eq6proof-theorem-global-adaptation}).

\item[$3^0$.]  \textit{Case} $(2\beta+1)p<r$. Set $q_k=r$ if
    $2^{-k}\le h^*$ and $q_k=p$ if $2^{-k}> h^*$, where
    $h^*$ will be chosen later. The following embedding holds [see
    page~62 in \citet{BesovIlNik}]:
    $\mathbb{N}_p(\beta,L)\subseteq\mathbb{N}_r (\beta -1/p+1/r,c_6L
    )$. Thus, by applying Lemma \ref{lemlp-norm-of-bais} with
    $\mathfrak{p}=r$, $s=\beta-1/p+1/r$ and $\mathcal{Q}=c_6L$, we
    obtain
%
\begin{eqnarray}
\label{eq7proof-theorem-global-adaptation}
\Xi_1&:=&\sum_{k\dvtx   2^{-k}\le h^*}
\biggl(\frac{
\ln(n)}{n2^{-k}} \biggr)^{(r-q_k)/2}\int \bigl( \Delta^*_{\mathcal{K},f}
\bigl(2^{1-k},y \bigr) \bigr)^{q_k}\,\mathrm{d}y
\nonumber\\[-8pt]\\[-8pt]
&=&\sum_{k\dvtx   2^{-k}\le h^*}\int \bigl( \Delta^*_{\mathcal{K},f}
\bigl(2^{1-k},y \bigr) \bigr)^{r}\,\mathrm{d}y \le
\bar{c}_9 L^{r}\bigl(h^*\bigr)^{\beta r - r/p +1}.\nonumber
\end{eqnarray}
By applying the same lemma 
with $\mathfrak{p}=r$, $s=\beta$ and $\mathcal{Q}=L$, we get
%
\begin{eqnarray}
\label{eq8proof-theorem-global-adaptation} \Xi_2&:=&\sum_{k\dvtx   2^{-k}> h^*}
\biggl(\frac{
\ln(n)}{n2^{-k}} \biggr)^{(r-q_k)/2}\int \bigl( \Delta^*_{\mathcal{K},f}
\bigl(2^{1-k},y \bigr) \bigr)^{q_k}\,\mathrm{d}y
\nonumber
\\
&=&\bar{c}_{10}L^p \bigl(n^{-1} \ln(n)
\bigr)^{(r-p)/2}\sum_{k\dvtx   2^{-k}> h^*}2^{-k
[\beta p- (r-p)/2 ]}
\\
&\le& \bar{c}_{11}L^p \bigl(n^{-1} \ln(n)
\bigr)^{(r-p)/2}\bigl(h^*\bigr)^{\beta p-(r-p)/2}.\nonumber
\end{eqnarray}
Here we have used that $\beta p-(r-p)/2<0$. In view of
(\ref{eq7proof-theorem-global-adaptation}) and
(\ref{eq8proof-theorem-global-adaptation}), we choose $h^*$ from the
equality $ L^{r}(h^*)^{\beta r -r/p +1}=L^p (n^{-1} \ln(n)
)^{(r-p)/2}(h^*)^{\beta p-(r-p)/2} $, so that $h^*= (L^{-2}n^{-1}\ln(n)
)^{1/(2\beta-2/p+1)}$. Finally, we obtain that
%
\begin{equation}
\label{eq9proof-theorem-global-adaptation} \qquad\quad \Xi\le\bar{c}_{12}L^{(r(1/2-1/r))/(\beta-1/p+1/2)}
\bigl(n^{-1}\ln(n) \bigr)^{(r(\beta-1/p+1/r))/(2\beta-2/p+1)}.
\end{equation}
The assertion of the theorem in the case $(2\beta+1)p<r$ follows now
from Theorem~\ref{thglobal-oracle-inequality},
(\ref{eq2proof-theorem-global-adaptation}) and
(\ref{eq9proof-theorem-global-adaptation}).
\end{longlist}\upqed
\end{pf*}

\section{Unknown design density}\label{secunknowndensity}
In this section we briefly comment on the case when the design density
$ g $ is unknown. We provide changes to be done in the selection rule
and in the presentation of the main result established in Theorems
\ref{thlocal-oracle-inequality} and \ref{thglobal-oracle-inequality}.
We also explain basic ideas related to the proofs of the new results.

In the context of the unknown design density, it is standard practice
to use a~plug-in estimator. This idea goes back to the Nadaraya--Watson
estimator and the problem considered in the paper is not an exception.

Suppose that an additional independent of an $ \{ X_i \}_{i=1}^n $
sample, say, $  \{ \widetilde{X}_i \}_{i=1}^n $, is available.
Alternatively, one can split the sample into two nonoverlapping parts.
Let $\widetilde{\mathbb{P}}_g^{(n)}$ stand for the probability law of $
\{ \widetilde{X}_i \}_{i=1}^n $. We reinforce
Assumption~\ref{assassumption-on-design} by the following condition:
$g\in\mathbb {G}\subset \{\ell\dvtx \mathbb{R}^2\to\mathbb{R}\dvtx
\|\ell\| _\infty\le \bar{g} \}$, where $\bar{g}<\infty$ and $\|\cdot\|
_\infty$ denotes the supremum norm on $\mathbb{R}^2$. It is noteworthy
that the constants $ \underline{g} $ and $ \bar{g} $ are both unknown.

Assume that based on $ \{ \widetilde{X}_i \}_{i=1}^n $ we can construct
an estimator of the design density $ g $, say, $ \hat g $, having the
following property. There exists a positive sequence $ a_n \downarrow0
$ as $ n \to\infty$ such that, for all sufficiently large $ n $,
%
\begin{equation}
\label{eq1density-est} \sup_{g\in\mathbb{G}}\widetilde{\mathbb{P}}_g^{(n)}
\bigl\{\| \hat g-g \|_\infty\ge a_n \bigr\}\le n^{-4r}.
\end{equation}
Denote by $ \Delta= [-3, 3]^2 $ the interval from
Assumption~\ref{assassumption-on-design} and introduce $
\underline{\hat g}=\inf_{x\in\Delta}\hat g(x) $ and $ \underline{\hat
g}{}_n=\underline{\hat g}\vee b_n $,
where $b_n$ tends to zero rather slowly. Theoretically, $b_n$ can be
chosen arbitrary, but a compromise allowing to keep our results under
reasonable sample size is $b_n=\ln^{-3}(n)$.

\subsection*{Changes in the selection rule (\protect\ref{eqselection-rule}) and in the oracle inequalities}
\begin{longlist}[${4^0}$.]
\item[${1^0}$.]  In the definition of estimators $\widehat F_{(\theta,
    h)} (\cdot)$ and $\widehat F_{(\theta, h) (\nu, h)}(\cdot)$, the
    unknown now values $g(X_i)$,  $ i=1,\ldots, n$, should be replaced by
    their truncated estimators $\hat g(X_i)\vee b_n$, $ i=1,\ldots,
    n$.

\item[${2^0}$.]  In the definition of all constants presented in
    Section~\ref{secconstants}, the unknown now value
    $\underline{g}^{-1}$ has to be replaced by $8\underline{\hat
    g}^{-2}\|\hat g\|_\infty$.

\item[${3^0}$.]  Let $\widehat{\operatorname{TH}}(\cdot)$ be obtained
    from $\operatorname{TH}(\cdot)$ by the replacement indicated
    in~${2^0}$. Then one should use
$ \operatorname{TH}^{(\mathrm{new})}(\eta)=\widehat{\operatorname
{TH}}(\eta)+2a_n \underline{\hat g}^{-1} \| \mathcal{K}\|_1^2 \widehat
F_{\infty}$ in (\ref{eqselection-rule}).

\item[${4^0}$.]  The right-hand sides of the local and global oracle
    inequalities established in Theorems
    \ref{thlocal-oracle-inequality} and
    \ref{thglobal-oracle-inequality} will additionally contain a term
    $c a_n$, where $c$ is a numerical constant independent of $F$, $g$
    and the sample size~$n$.
\end{longlist}

\subsection*{Sketch of the proof of the new version of Theorem~\protect\ref{thlocal-oracle-inequality}}
\begin{longlist}[3.]
\item[1.]  Denote\vspace*{1pt} by $ \mathcal C$ the event $\{\|\hat g-g \|_\infty\le
    a_n
    \} $. Similar to the proof given in the step \textit{Risk
    computation under $ \widebar{\mathcal B} $} of Theorem
    \ref{thlocal-oracle-inequality}, the computations under the event $
    \widebar{\mathcal C} $ lead to the same reminder term in the oracle
    inequality.

\item[2.]  All computations under the event $ \mathcal C $ are done
    conditionally with respect to $ \{ \widetilde{X}_i \}_{i=1}^n $. It
    allows us to treat the estimator $ b_n\vee\hat g(\cdot) $ as
    nonrandom.
\begin{longlist}[2.2.]
\item[2.1.]  Analyzing the proof of ``probabilistic'' Lemmas
    \ref{lemgauss-on-matrices} and \ref{leminegality-for-sup-norm}, we
    see that the results remain valid if we replace the function $ g $
    in the denominator of all expressions by an arbitrary function
    bounded from below and above on~$ \Delta$. Thus, the use of $ b_n
    \vee\hat g (\cdot) $ in place of $ g (\cdot) $ under $ \mathcal
    C $ is eligible. That leads to the similar assertion where $
    \underline{g}^{-1/2} $ is substituted with $ \underline{\hat
    g}{}_n^{-1}\|g\|_\infty^{-1/2} $ in all the constants\vspace*{-2pt} involved. The
    latter quantity can be bounded under $ \mathcal C $ by $ 2
    \underline{\hat g}^{-1}\|g\|_\infty^{-1/2} $, for $ n $ large
    enough.

However, such quantities cannot be used in the definition of the
threshold directly, because they incorporate the unknown $ \|g\|_\infty
$. Nevertheless, under the event~$ \mathcal C $,  $\|g\|_\infty$ can be
bounded by $ 2\|\hat g\|_\infty$. Moreover, we remark that all the
quantities listed in Section~\ref{secconstants} are increasing
functions of $\underline{g}^{-1}$ so we can\vspace*{-2pt} replace $
\underline{g}^{-1} $ by the upper bound $ 8 \underline{\hat
g}^{-2}\|\hat g\|_\infty$ available under $ \mathcal C $. It
explains~${2^0}$.

\item[2.2.]  The replacement of $ g $ by $ b_n \vee\hat g $
    leads
    to an additional ``approximation error'' bounded from above by $
    A(g):=\sup_{x\in\Delta} |g(x) [b_n\vee\hat g (x)]^{-1} - 1 |.
$ This quantity should be added to $ \widehat{\operatorname{TH}}(\cdot
) $ in order to preserve the proof of
Theorem~\ref{thlocal-oracle-inequality}. Since $ A(g) $ depends on $g$,
one should instead use its upper bound $ a_n \underline{\hat g}^{-1} $
which is available under $ \mathcal C $.
It gives ${3^0}$ for $n$ large enough.

\item[2.3.]  The use of $\operatorname{TH}^{(\mathrm{new})}(\cdot)$
    in place of $\operatorname{TH} (\cdot)$ leads to an additional term
    $ 2a_n \underline{\hat g}^{-1} \| \mathcal{K}\|_1^2  (3M +
    4C_5^{(\mathrm{new})} )$ in (\ref{eq020proof-of-theorem-local}). It
    explains ${4^0}$ and completes the sketch of the proof of Theorem
    \ref{thlocal-oracle-inequality}. Since the global oracle inequality
    is obtained by the integration of the local one over bounded
    interval of $\mathbb{R}^2$, the assertion of Theorem
    \ref{thglobal-oracle-inequality} remains the same up to the term
    $ca_n$.
\end{longlist}
\end{longlist}

\subsection*{The additional assumption about $g$ and an example of an estimator obeying~(\protect\ref{eq1density-est})}

If we suppose that $\mathbb{G}\subseteq\mathbb{H}_2(\gamma,R)$, where
$\mathbb{H}_2(\gamma,R)$ is an isotropic H\"older class of two-variate
functions, then $ a_n= [ n^{-1} \ln(n)  ]^{\gamma/(2(\gamma+1))}$,  and
this rate is attainable by a kernel estimator with properly chosen
kernel and bandwidth. This yields together with ${4^0}$ that if $
\gamma>2 \mathbf{b}$, the adaptive results established in Theorems
\ref{thpointwise-adaptation} and \ref{thglobal-adaptation} remain
unchangeable.

Another possibility is to suppose that $\mathbb{G}$ is a parametric
family of densities. In this case, $c a_n$ can be viewed as a reminder
term.

\section*{Acknowledgments}
The authors wish to thank an Associate Editor and anonymous referees
for their comments that substantially improved the paper.

\begin{supplement}
\stitle{Proofs of lemmas for
``Adaptive estimation under single-index constraint in a
regression model''}
\slink[doi]{10.1214/13-AOS1152SUPP}
\sdatatype{.pdf}
\sfilename{aos1152\_supp.pdf}
\sdescription{We provide detailed proofs of the auxiliary results (Lemmas \ref{lemboundsforbias}--\ref{leminegality-for-sup-norm} and
\ref{lemlp-norm-of-bais}) for the paper.}
\end{supplement}





\printaddresses

\end{document}